\documentclass[]{article}
\usepackage{graphicx}
\usepackage{amsfonts}
\usepackage{amsmath,amssymb,amscd,mathrsfs,amsthm}
\usepackage{fancyhdr}
\usepackage{mathrsfs}
\usepackage{float}
\usepackage{url}
\usepackage{titlesec}
\usepackage{indentfirst}
\usepackage{booktabs}
\usepackage{verbatim}
\usepackage{color}
\usepackage[colorlinks,linkcolor=black,anchorcolor=black,citecolor=black]{hyperref} 
\usepackage{pdflscape}

\usepackage{todonotes}

\usepackage[page,header]{appendix}
\usepackage{titletoc}
\usepackage[noabbrev]{cleveref}
\usepackage{bbm}
\newcommand{\rd}{\,\mathrm{d}}

\newcommand{\abs}[1]{\lvert #1 \rvert}

\newcommand{\p}{\partial}

\newcommand{\fl}[2]{\frac{#1}{#2}}

\newcommand{\dt}{\Delta}
\newcommand{\suml}{\sum\limits}

\newcommand{\be}{\begin{equation}}
\newcommand{\ee}{\end{equation}}
\newcommand{\bal}{\begin{aligned}}
\newcommand{\eal}{\end{aligned}}
\newcommand{\ba}{\begin{array}}
\newcommand{\ea}{\end{array}}
\newcommand{\bea}{\begin{eqnarray}}
\newcommand{\eea}{\end{eqnarray}}
\newcommand{\beas}{\begin{eqnarray*}}
\newcommand{\eeas}{\end{eqnarray*}}
\newcommand{\beit}{\begin{itemize}}
\newcommand{\eit}{\end{itemize}}

\newcommand{\bmat}[1]{\begin{bmatrix} #1 \end{bmatrix}}

\renewcommand{\vec}[1]{\boldsymbol{\mathrm{#1}}}

\providecommand{\vn}{\ensuremath{\vec{n}}}

\providecommand{\vu}{\ensuremath{\vec{u}}}
\providecommand{\vv}{\ensuremath{\vec{v}}}

\providecommand{\vx}{\ensuremath{\vec{x}}}

\numberwithin{equation}{section}

\begin{document}
	\title{An adaptive dynamical low rank method for the nonlinear Boltzmann equation\footnote{This work is partially supported under the NSF CAREER grant DMS-2153208, NSF CDS\&E grant CBET-1854829, and AFOSR grant FA9550-21-1-0358.}}
	
    \author{
	    Jingwei Hu\footnote{Department of Applied Mathematics, University of Washington, Seattle, WA 98195, USA (hujw@uw.edu).} \  \
	    and \ Yubo Wang\footnote{Department of Mathematics, Purdue University, West Lafayette, IN 47907, USA (wang3158@purdue.edu).}}  
	\maketitle

\begin{abstract}


Efficient and accurate numerical approximation of the full Boltzmann equation has been a longstanding challenging problem in kinetic theory. This is mainly due to the high dimensionality of the problem and the complicated collision operator. In this work, we propose a highly efficient adaptive low rank method for the Boltzmann equation, concerning in particular the steady state computation. This method employs the fast Fourier spectral method (for the collision operator) and the dynamical low rank method to obtain computational efficiency. An adaptive strategy is introduced to incorporate the boundary information and control the computational rank in an appropriate way. Using a series of benchmark tests in 1D and 2D, we demonstrate the efficiency and accuracy of the proposed method in comparison to the full tensor grid approach.

\end{abstract}




\section{Introduction}
Kinetic theory describes the non-equilibrium dynamics of gases or systems comprised of a large number of particles. It provides rich information at the mesoscopic level when the well-known fluid mechanical laws of Navier-Stokes and Fourier become inadequate. Various applications of kinetic theory can be found in fields such as rarefied gas dynamics \cite{cercignani2000rarefied}, plasma physics \cite{birdsall2004plasma}, semiconductor modeling \cite{markowich2012semiconductor} and biological and social sciences \cite{naldi2010mathematical}. 

In this work, we are interested in the efficient numerical approximation of the nonlinear Boltzmann equation \cite{Cercignani,Villani02}, which is the central model in kinetic theory and reads as
\be
\p_t f+\vv\cdot \nabla_{\vx} f =\mathcal{Q}(f,f),\quad  t>0, \quad \vx\in \Omega_{\vx}\subset\mathbb{R}^d, \quad \vv\in\mathbb{R}^d,
\label{eqn:Boltzmann}
\ee
where $f = f(\vx,\vv,t)$ is the phase space distribution function of time $t$, position $\vx$, and particle velocity $\vv$; $\mathcal{Q}$ is the Boltzmann collision operator, which is a quadratic integral operator modeling the binary interactions between particles. It is convenient to introduce the bilinear form of $\mathcal{Q}$:
\be
\bal
\mathcal{Q}(g, f)(\vv) &= \int_{\mathbb{R}^d}\int_{S^{d-1}} B(|\vv-\vv_*|,\cos \chi) [g(\vv_{*}')f(\vv')-g(\vv_*)f(\vv)]\, {\rm d} \sigma \,  {\rm d} \vv_*,
\label{eqn:collision_operator}
\eal 
\ee
where the post-collisional velocities $(\vv',\vv_*')$ are defined in terms of pre-collisional velocities $(\vv,\vv_*)$ through the conservation of momentum and energy during the collision:
\be
\vv^{\prime}=\frac{\vv+\vv_{*}}{2}+\frac{\left|\vv-\vv_{*}\right|}{2} \sigma, \quad \vv_{*}^{\prime}=\frac{\vv+\vv_{*}}{2}-\frac{\left|\vv-\vv_{*}\right|}{2} \sigma,
\label{eqn:collision velocity}
\ee
with $\sigma$ being a vector over the unit sphere $S^{d-1}$. The collision kernel $B$ is a non-negative function depending on $|\vv-\vv_*|$ and cosine of the deviation angle $\chi$, i.e., the angle between $\vv-\vv_*$ and $\vv'-\vv_*'$. It should be noted that collisions happen only in the velocity space, thus time and spatial dependence is omitted in (\ref{eqn:collision_operator}). This property has important consequence in design of efficient numerical methods as we shall see later. With the distribution function $f$, one can retrieve the macroscopic quantities via its moments:
\be
\int_{\mathbb{R}^d} f(\vx,\vv,t)\bmat{1 \\ \vv \\ \fl {1} 2 |\vv|^2}{\rm d}\vv=\bmat{\rho(\vx,t)\\ \rho(\vx,t) {\bf u}(\vx,t) \\ \fl 1 2 {\rho(\vx,t)|{\bf u}(t,\vx)|^2+\fl {{d}}{2}\rho(\vx,t) RT(\vx,t)}},
\ee
where $\rho(\vx,t)$, ${\bf u}(\vx,t)$, and $T(\vx,t)$ are the density, bulk velocity, and temperature, $R$ is the Boltzmann constant.

Despite of the long history and wide application of the Boltzmann equation, numerically solving the Boltzmann equation still faces great challenges nowadays. This is mainly due to the high dimensionality of the equation and the complicated collision operator. The prevailing method is the direct simulation Monte Carlo (DSMC) method \cite{nanbu1980direct,bird1994molecular} because it can avoid the curse of dimensionality. DSMC method models binary collisions stochastically but could suffer from slow convergence in certain cases such as low speed or near continuum flows. On the other hand, the deterministic method based on discretization of the equation on representative grids has undergone significant development over the past decade. This is partly due to the rapid growth of the computing power as well as the algorithmic advance in approximation of the Boltzmann collision operator. Regarding the latter, the Fourier spectral method \cite{PP96, PR00} stands out for its high accuracy and possibility of being further accelerated by the fast Fourier transform (FFT). The readers can refer to \cite{Pareschi, Hu21} for a review of such methods. Relevant to the current work, we mention the fast algorithm proposed in \cite{MP06} which can efficiently evaluate the collision operator for certain collision kernel in $\mathcal{O}(MN_{\vv}^d\log N_{\vv})$ complexity, where $N_{\vv}$ is the number of points in each velocity dimension and $M$ is the number of points over $S^{d-1}$. Even equipped with the fast solver for the collision operator, solving the Boltzmann equation deterministically can still be very expensive. In the full tensor grid approach, the overall complexity (per time step) would be $\mathcal{O}(N_{\vx}^dMN_{\vv}^{d}\log N_{\vv})$, where $N_{\vx}$ is the number of discretization points used in each spatial dimension. This motivates us to seek more efficient method to overcome the intrinsic high dimensionality of the problem.

Recently, a class of dynamical low rank method has been applied to solving kinetic equations including the Vlasov equation \cite{einkemmer2018low,einkemmer2019quasi}, Boltzmann-BGK equation \cite{einkemmer2019low, einkemmer2021efficient} and radiation transfer equation \cite{einkemmer2021asymptotic,peng2020low}. The basic idea is to find a low-rank approximation of the unknown function $f$ by projecting the equation onto the tangent space of the low-rank solution manifold. Upon a further operator splitting, the original $N_{\vx}^dN_{\vv}^d$ dimensional problem can be reduced to a few $N_{\vx}^d$ or $N_{\vv}^d$ dimensional problems. We mention that this formulation can only be done easily if the original equation has a tensor structure ($\vx$ and $\vv$ are separated in some sense). For some collision operators, e.g., the BGK operator which is highly nonlinear due to the $f$ dependent Maxwellian function, extra effort is needed to make the method efficient \cite{einkemmer2021efficient}.

In this work, we apply the dynamical low rank method to the Boltzmann equation (\ref{eqn:Boltzmann}) and investigate the performance of the method in a series of benchmark tests that concern the steady state solutions. Our contribution can be summarized as follows: 1) The dynamical low rank method is for the first time applied to the nonlinear Boltzmann equation (\ref{eqn:Boltzmann}). Even though the collision operator $\mathcal{Q}(f,f)$ is quite complicated, it is local in $\vx$ hence is highly suited in the low rank framework. Moreover, the previously developed fast Fourier spectral method can be applied straightforwardly to accelerate the overall method. 2) We propose an adaptive strategy to add and remove basis along the time evolution. In particular, the stage of adding the basis is strongly motivated by our underlying problem. Most benchmark tests for the Boltzmann equation involve steady state solutions (e.g., normal shock, Couette flow, thermally driven cavity flow, etc. \cite{JAH19_1}) for which the boundary condition is highly non-trivial and plays an important role. We show that to accurately simulate this type of problems, the boundary information needs to be added to the solution on the fly. As a consequence, dropping the basis becomes mandatory, otherwise the numerical rank will increase constantly. This is in contrast to most of the previous dynamical low rank methods on kinetic equations, where a fixed rank can often be used throughout the simulation. 3) Using asymptotic analysis and heuristic arguments, we identify a class of problems -- normal shock problem -- whose steady state solutions are indeed low rank in some regimes, and further confirm it in numerical experiments. This provides some theoretical guarantee for the proposed low rank method to be an efficient approach for solving the nonlinear Boltzmann equation.

The rest of this paper is organized as follows. In Section \ref{sec:dynamical}, we describe the dynamical low rank method for the Boltzmann equation, including the time, velocity and physical space discretization as well as the treatment of the boundary condition. In Section \ref{sec:adaptive}, we introduce an adaptive strategy to add and drop basis in the dynamical low rank method during the time evolution. In Section \ref{sec:lowrank}, we analyze the normal shock problem and demonstrate the low rank property of the solution in both the weak and strong shock wave regimes. Section \ref{sec:num} presents numerical examples in 1D and 2D using the proposed adaptive dynamical low rank method. Several benchmark tests for the nonlinear Boltzmann equation are considered: normal shock, Fourier flow, lid driven cavity flow, and thermally driven cavity flow. The paper is concluded in Section \ref{sec:con}.

\section{The dynamical low rank method for the Boltzmann equation}
\label{sec:dynamical}

In this section, we introduce the dynamical low rank method for the Boltzmann equation (\ref{eqn:Boltzmann}). We first present the formulation in the continuous setup, where we highlight the special structure of the collision operator in obtaining an efficient low rank approximation. We then describe the discretization in the velocity space and physical space, and treatment of the typical boundary conditions of the Boltzmann equation. Finally, we add the time discretization to obtain a fully discrete low rank scheme.

The starting point of the method is to constrain the distribution function $f(\vx,\vv,t)$ to a low rank manifold $\mathbb{M}$ such that 
\be \label{lowrank decomposition:stat}
	f(\vx,\vv,t) = \sum_{i,j=1}^rX_i(\vx,t)S_{ij}(t)V_j(\vv,t),
\ee
where $r$ is the representation rank and the basis functions $\{X_i\}_{1\leq i\leq r} \subset L^2(\Omega_{\vx})$ and $\{V_j\}_{1\leq j \leq r} \subset L^2(\Omega_{\vv})$ are orthonormal:
\be
\label{eqn:orthogonality}
{\langle X_i,X_j\rangle}_{\vx}=\delta_{ij}, \quad  {\langle V_i,V_j\rangle}_{\vv}=\delta_{ij},\quad 1 \le i,j \le r,
\ee
with $\langle \cdot,\cdot\rangle_{\vx}$ and $\langle \cdot,\cdot\rangle_{\vv}$ being the inner products on $L^2(\Omega_{\vx})$ and $L^2(\Omega_{\vv})$, respectively. Note here we consider a finite velocity domain $\Omega_{\vv}$ rather than the whole space $\mathbb{R}^d$ to avoid the complication in the infinite domain. This is a reasonable assumption because the majority of the numerical methods for kinetic equations need to first truncate the velocity domain and then perform the discretization. It can often be done without much loss of accuracy since $f$ decays sufficiently fast as $\vv$ goes to infinity.

We rewrite  \cref{eqn:Boltzmann} as
\be \label{eqn:RHS}
\p_tf=-\vv \cdot \nabla_{\vx} f+\mathcal{Q}(f,f) :=\text{RHS}.
\ee
To ensure uniqueness of the dynamical factors $X_i$, $S_{ij}$, and $V_j$ through \cref{eqn:RHS}, we impose the following gauge conditions by constraining the derivatives in the null space (for details, see \cite{koch2007dynamical}):
\be \label{eq:gauge}
{\langle \partial_t X_i,X_j\rangle}_{\vx}=0, \quad {\langle \partial_t V_i,V_j\rangle}_{\vv}=0,\quad 1 \le i,j \le r.
\ee
We now project the right hand side of  (\ref{eqn:RHS}) onto the tangent space of $\mathbb{M}$:
\be
	\p_t f = P_f(\text{RHS}),
\label{projected equation}
\ee
where the orthogonal projector $P_f$ can be written as
\be
	P_f(\text{RHS}) = \sum_{j=1}^r\langle V_j,\text{RHS}\rangle_{\vv} V_j - \sum_{i,j=1}^rX_i\langle X_iV_j,\text{RHS}\rangle_{\vx,\vv}V_j+\sum_{i=1}^rX_i\langle X_i,\text{RHS}\rangle_{\vx}. \label{projected_parts}
\ee

To avoid the possible ill-conditioning of the matrix $S = (S_{ij})_{1 \le i,j \le r}$, one can perform a simple operator splitting \cite{lubich2014projector} to decompose (\ref{projected_parts}) into three subflows:
\begin{align}
    \partial_t f &= \sum_{j=1}^r\langle V_j,\text{RHS}\rangle_{\vv} V_j, \label{eqn:K-step general}\\
    \partial_t f &= -\sum_{i,j=1}^rX_i\langle X_iV_j,\text{RHS}\rangle_{\vx,\vv}V_j, \label{eqn:S-step general}\\
    \partial_t f &= \sum_{i=1}^rX_i\langle X_i,\text{RHS}\rangle_{\vx}.\label{eqn:L-step general}
\end{align}
Using the orthogonality condition (\ref{eqn:orthogonality}) and the gauge condition (\ref{eq:gauge}), we can further simplify each subflow and proceed in the following three substeps:
\begin{itemize}
\item $K$-step: Define $K_{j}(\vx,t)=\suml_{i=1}^{r} X_{i}(\vx,t) S_{i j}(t)$
then $f(\mathbf{x}, \mathbf{v},t)=\suml_{j=1}^{r} K_{j}(\vx,t) V_{j}(\vv,t).$ We can rewrite \cref{eqn:K-step general} as
\begin{equation} \label{eqn:KK}
\partial_t \left(\sum_{j=1}^rK_j V_j\right)=\sum_{j=1}^r \left(\partial_t K_j V_j +K_j \partial_tV_j\right)=  \sum_{j=1}^r\langle V_j,\text{RHS}\rangle_{\vv} V_j.
\end{equation}
Using the orthogonality of $\{V_j\}_{1\le j \le r}$ and ${\langle \partial_t V_j,V_k\rangle}_{\vv}=0$ for $1 \le j,k \le r$, we have
\begin{equation}
\begin{aligned}
\partial_{t} K_{j} &=\left\langle V_{j}, \mathrm{RHS}\right\rangle_{\mathbf{v}}  \\
&=-\sum_{l=1}^{r}\left\langle \vv V_{j} V_{l}\right\rangle_{\mathbf{v}} \cdot \nabla_{\mathbf{x}} K_{l}+\sum_{m, n=1}^{r}\left\langle V_{j} \mathcal{Q}\left(V_{m}, V_{n}\right)\right\rangle_{\mathbf{v}} K_{m} K_{n}, \quad j =1,\dots,r,\label{eqn:K-step detailed}
\end{aligned}
\end{equation}
where the simplification of the last term relies crucially on the bilinearity of the collision operator (\ref{eqn:collision_operator}) as well as the fact that collisions act locally in the physical space. It can be seen that (\ref{eqn:K-step detailed}) together with $\partial_tV_j=0$ solve (\ref{eqn:KK}). Since the solution to the subflow is unique, we thus know $\{V_j\}_{1\le j \le r}$ remains unchanged during this substep. 
\item $S$-step: We can argue similarly to obtain that the subflow (\ref{eqn:S-step general}) is equivalent to
\begin{equation}
\begin{aligned}
\partial_t S_{ij} &= -\langle X_iV_j,\text{RHS}\rangle_{\vx,\vv}\\
&=\sum\limits_{k,l=1}^r\langle \vv V_j V_l\rangle_{\vv}\cdot \langle X_i \nabla_{\vx}X_k\rangle_{\vx}S_{kl}-\sum\limits_{k,l,m,n=1}^r\langle X_i X_kX_l\rangle_{\vx}\langle V_j \mathcal{Q}(V_m,V_n)\rangle_{\vv}S_{km}S_{ln}, \quad i,j=1,\dots,r.\label{eqn:S-step detailed}
\end{aligned}
\end{equation}
During this substep, both $\{V_j\}_{1\le j \le r}$ and $\{X_i\}_{1\le i \le r}$ remain unchanged.
\item $L$-step: Define $L_{i}(\vv,t)=\suml_{j=1}^{r} S_{i j}(t) V_{j}(\vv,t)$
then $f(\mathbf{x}, \mathbf{v},t)=\suml_{i=1}^{r} X_{i}(\mathbf{x},t) L_{i}(\mathbf{v},t) .$ By similar arguments, the subflow (\ref{eqn:L-step general}) is equivalent to
\begin{equation}
\begin{aligned}
\partial_t L_i  &= \langle X_i,\text{RHS}\rangle_{\vx},  \\
&=-\sum\limits_{l=1}^r \vv\cdot \langle X_i \nabla_{\vx} X_l\rangle_{\vx} L_l+\sum\limits_{m,n = 1}^r\mathcal{Q}(L_m,L_n)\langle X_i X_mX_n\rangle_{\vx}, \quad i=1,\dots,r.
\label{eqn:L-step detailed}
\end{aligned}
\end{equation}
During this substep, $\{X_i\}_{1\le i \le r}$ remains unchanged.
\end{itemize}

Therefore, we have obtained a set of low rank equations (\ref{eqn:K-step detailed})-(\ref{eqn:L-step detailed}) in the continuous setting. The task remains is to apply the proper discretization to these equations in the velocity space, physical space, and time, which we will detail in the following subsections.

\subsection{Velocity space discretization} 
\label{velocity discretization}

Examining the equations (\ref{eqn:K-step detailed})-(\ref{eqn:L-step detailed}), we can see that all terms pertaining to the collision operator have the form of $\mathcal{Q}(h_1,h_2)$, where $h_1$ and $h_2$ are some functions of $\vv$. Luckily this isn't much change from the original collision operator in (\ref{eqn:Boltzmann}) and we can apply the well-developed fast Fourier spectral methods. 

Specifically, for 2D Maxwell molecules ($d=2$ and $B=\text{const}$) and 3D hard spheres ($d=3$ and $B=\text{const}|\vv-\vv_*|$), we can use the algorithm proposed in \cite{mouhot2006fast} with complexity $\mathcal{O}(MN_{\vv}^{d}\log N_{\vv})$, where $N_{\vv}$ is the number of points in each dimension of the velocity space and $M\ll N_{\vv}^{d-1}$ is the number of points used on the sphere $S^{d-1}$; for general collision kernel $B=B(|\vv-\vv_*|,\cos \chi)$, we can use the algorithm proposed in \cite{gamba2017fast} with complexity $\mathcal{O}(MN_{\vv}^{d+1}\log N_{\vv})$. Both algorithms can be implemented as a discrete velocity method: one chooses an appropriate velocity domain $[-L_{\vv},L_{\vv}]^d$ and uniform grid points $\{\vv_q\}$; the collision solver takes discrete values $\{h_1(\vv_q)\}$ and $\{h_2(\vv_q)\}$ and outputs $\{\mathcal{Q}(h_1,h_2)(\vv_q)\}$ on the same set of grid points. For more details, the readers can refer to \cite{mouhot2006fast, gamba2017fast}.


\subsection{Physical space discretization}
\label{space discretization}

There are various ways to discretize the equations (\ref{eqn:K-step detailed})-(\ref{eqn:L-step detailed}) in the physical space, for example, one can apply the Fourier spectral method \cite{einkemmer2018low} or the high resolution finite difference scheme \cite{einkemmer2021efficient} directly to these equations. Generally speaking, the conventional scheme used for the original equation needs to be tailored when solving the equations resulted from the low rank projection. The issue also becomes a bit tricky when the boundary condition is not periodic. 

Here we adopt a ``first discretize, then project" strategy, which is simpler because it follows directly from the scheme for the original equation. We mention that this idea is similar to the so-called kinetic flux vector splitting (KFVS) scheme \cite{Deshpande86}, a well-known method for solving the compressible Euler equations derived from the kinetic equation. For simplicity, we focus on the first order upwind scheme in this work. To extend it to high order, similar strategy for the KFVS scheme \cite{mandal1994kinetic} can be considered. 

We use the one-dimensional case ($d=1$) to illustrate the idea. Extension to high dimension with rectangular grid is straightforward as implemented in our numerical examples. Assume $\Omega_{x} = [-L_x,L_x]$ with uniform grid points chosen as $x_p = -L_x+ (p-\fl 1 2)\dt x$, $p= 1,\ldots,N_{x}$, $\Delta x= \frac{2L_x}{N_{x}}$. Since the transport term in the Boltzmann equation (\ref{eqn:Boltzmann}) is linear, it is very easy to apply the upwind scheme:
\be
\bal
\p_t f(x,v,t)=& -\fl {v+\abs{v}}{2}\fl {f(x,v,t)-f(x-\Delta x,v,t)}{\Delta x}\\
&-\fl {v-\abs{v}}{2}\fl {f(x+\Delta x,v,t)-f(x,v,t)}{\Delta x} + \mathcal{Q}(f(x,v,t),f(x,v,t))\\
:=&-v^+ D_+f(x,v,t)-v^-D_-f(x
,v,t) + \mathcal{Q}(f(x,v,t),f(x,v,t)),
\label{scheme:fullupwind}
\eal
\ee
where $v^{\pm}=\fl {v \pm \abs{v}}{2}$, and $D_{\pm}$ are first order upwind operators. 

For \cref{scheme:fullupwind}, we can apply the same projection process as we did previously to \cref{eqn:RHS} to obtain (i.e., the analogs of (\ref{eqn:K-step detailed})-(\ref{eqn:L-step detailed})):
\begin{itemize}
	\item $K$-step: 
	\begin{equation}
	\bal
	\p_t K_j(x,t) =& -\sum_{l = 1}^r\langle v^+ V_j(v,t)V_l(v,t)\rangle_{v}D_+K_l(x,t) -\sum_{l= 1}^r\langle v^- V_j(v,t)V_l(v,t)\rangle_{v}D_-K_l(x,t)\\
	&+\sum\limits_{m,n = 1}^r\langle V_j(v,t)\mathcal{Q}(V_m(v,t),V_n(v,t))\rangle_{v} K_m(x,t)K_n(x,t). 
	\eal
	\end{equation}
	\item $S$-step: 
        \begin{equation}
	\bal
	\partial_t S_{ij}(t) = & \sum_{k,l= 1}^r\langle v^+ V_j(v,t)V_l(v,t)\rangle_{v}\langle X_i(x,t) D_+X_k(x,t)\rangle_{x}S_{kl}\\
	&+\sum_{k,l = 1}^r\langle v^- V_j(v,t)V_l(v,t)\rangle_{v}\langle X_i(x,t)D_-X_k(x,t)\rangle_{x}S_{kl}\\
	&-\sum\limits_{k,l,m,n= 1}^r\langle X_i(x,t)X_k(x,t)X_l(x,t)\rangle_{x}\langle V_j(v,t) \mathcal{Q}(V_m(v,t),V_n(v,t))\rangle_{v}S_{km}S_{ln}.
	\eal
	\end{equation}
	\item $L$-step:
	\begin{equation}
	\bal
	\p_t L_i(v,t)=&-\sum\limits_{l=1}^r v^+\langle X_i(x,t) D_+ X_l(x,t)\rangle_x L_l(v,t)-\sum\limits_{l=1}^r v^-\langle X_i(x,t)D_- X_l(x,t)\rangle_x L_l(v,t)\\
	&+\sum\limits_{m,n = 1}^r\mathcal{Q}(L_m(v,t),L_n(v,t))\langle X_i(x,t)X_m(x,t)X_n(x,t)\rangle_x.
	\eal
        \end{equation}
	\end{itemize}

\subsection{Treatment of the boundary condition}
\label{sec:boundary}
In the low rank framework, boundary condition for $f(\vx,\vv,t)$ needs to be transformed to the boundary condition of $\{K_j\}_{1\le j \le r}$. In fact, this transformation has a non-trivial impact on the fully discrete scheme which we shall describe in the next subsection.

For a boundary point $\vx \in \p\Omega_{\vx}$ with outward pointing normal $\vn(\vx)$ and boundary velocity $\vu_{w}(\vx,t)$, general boundary conditions for Boltzmann equation (\ref{eqn:Boltzmann}) are defined through the inflow direction:
\be
f(\vx,\vv,t) = f_{bdy}(\vx,\vv,t), \quad (\vv-\vu_w(\vx,t))\cdot \vn(\vx) < 0,
\ee
where $f_{bdy}$ is a prescribed function. The other half of $f(\vx,\vv,t)$ is given from interior of the domain (outflow). We thus define
\begin{equation}
f^b(\vx,\vv,t) = 
\left\{
\begin{aligned}
&f_{bdy}(\vx,\vv,t),& \quad &(\vv-\vu_w(\vx,t))\cdot \vn(\vx) < 0,\\
&f(\vx,\vv,t),&\quad &(\vv-\vu_w(\vx,t))\cdot \vn(\vx) \ge 0.
\end{aligned}
\right.
\label{eqn:full bdy general}
\end{equation}
Accordingly, we can project the full boundary $f^b(\vx,\vv,t)$ to the space spanned by $\{V_j\}_{1 \le j \le r}$ to obtain boundary values for $\{K_j\}_{1\le j \le r}$:
\begin{equation}
\begin{aligned} \label{boundary general}
K_j(\vx,t)=& \langle f^b(\vx,\vv,t),V_j(\vv,t)\rangle_{\vv}\\
=& \langle f_{bdy}(\vx,\vv,t)\mathbbm{1}_{(\vv-\vu_w(\vx,t))\cdot \vn(\vx) < 0}, V_j(\vv,t)\rangle_{\vv}+\langle f(\vx,\vv,t)\mathbbm{1}_{(\vv-\vu_w(\vx,t))\cdot \vn(\vx) \ge 0}, V_j(\vv,t)\rangle_{\vv}\\
=&\langle f_{bdy}(\vx,\vv,t)\mathbbm{1}_{(\vv-\vu_w(\vx,t))\cdot \vn(\vx) < 0}, V_j(\vv,t)\rangle_{\vv}+\sum_{l=1}^r K_l(\vx,t)\langle \mathbbm{1}_{(\vv-\vu_w(\vx,t)) \cdot \vn(\vx)\geq 0}V_l(\vv,t)V_j(\vv,t)\rangle_{\vv},
\end{aligned}
\end{equation}
where the $K_l$ term appearing on the right hand side of (\ref{boundary general}) can be approximated using values inside the domain (extrapolation) since the term results from the outflow.

Two typical boundary conditions used when solving the Boltzmann equation (\ref{eqn:Boltzmann}) are the following inflow boundary and Maxwell diffusive boundary. For inflow boundary, we take $\vu_w(\vx,t) = \textbf{0}$ and 
\begin{equation}
f_{bdy}(\vx,\vv,t) =  \frac{\rho_{in}(\vx,t)}{(2\pi T_{in}(\vx,t))^{d/2}}\exp\left(-\frac{|\vv-\vu_{in}(\vx,t)|^2}{2RT_{in}(\vx,t)}\right), \quad \vv \cdot \vn(\vx) < 0,
\end{equation}
where $\rho_{in}$, $\vu_{in}$ and $T_{in}$ are the density, bulk velocity and temperature of the prescribed inflow.
For the Maxwell diffusive boundary, we take
\begin{equation}
\begin{aligned}
f_{bdy}(\vx,\vv,t) &=  \rho_w(\vx,t)\exp\left(-\frac{|\vv-\vu_w(\vx,t)|^2}{2RT_w(\vx,t)}\right), \quad (\vv-\vu_w(\vx,t)) \cdot \vn(\vx) < 0,
\end{aligned}
\end{equation}
where $T_w$ is the wall temperature, and $\rho_w$ is determined by conservation of mass through the wall:
\begin{equation}
\rho_w(\vx,t) = -\frac{\int_{(\vv-\vu_w(\vx,t)) \cdot \vn(\vx) \ge 0} (\vv-\vu_w(\vx,t)) \cdot \vn(\vx) f(\vx,\vv,t)\, {\rm d} \vv}{\int_{(\vv-\vu_w(\vx,t)) \cdot \vn(\vx) < 0}(\vv-\vu_w(\vx,t)) \cdot \vn(\vx) \exp\left(-\frac{|\vv-\vu_w(\vx,t)|^2}{2RT_w(\vx,t)}\right)\, {\rm d} \vv}.
\end{equation}

\subsection{Time discretization and the fully discrete scheme}
\label{scheme:first order}

We now add the time discretization to (\ref{eqn:K-step detailed})-(\ref{eqn:L-step detailed}) to obtain a fully discrete scheme. Since most of the examples we are interested in this paper concern the stationary Boltzmann equation, the first order time discretization suffices. For high order method in time, the readers can refer to \cite{einkemmer2021asymptotic} and references therein.

Given the initial condition $f(\vx,\vv,0)=f^0(\vx,\vv)$, we first perform the singular value decomposition $f^0(\vx,\vv)=\sum_{i,j=1}^r X_i^0(\vx)S_{ij}^0V^0_j(\vv)
$ to obtain $(X_i^0,S_{ij}^0,V_j^0)$, where a fixed, reasonable rank $r$ is chosen and used in the following computation.

Suppose at time step $t^n$, $(X_i^n,S_{ij}^n,V_j^n)$ are available. In order to obtain $(X_i^{n+1},S_{ij}^{n+1},V_j^{n+1})$ at $t^{n+1}$, we proceed as follows:
\begin{enumerate}
\item {\bf $K$-step}. 
\begin{enumerate}
\item Construct $K^n_j = \sum_{i = 1}^r X_i^nS_{ij}^n$.
\item Perform the forward Euler step in (\ref{eqn:K-step detailed}) to obtain $K_{j}^{n+1}$:
\begin{equation}
\bal
K_j^{n+1} &=K_j^n-\Delta t\sum_{l=1}^{r}\left\langle \vv V_{j}^n V_{l}^n\right\rangle_{\mathbf{v}} \cdot \nabla_{\mathbf{x}} K_{l}^n+\Delta t \sum_{m, n=1}^{r}\left\langle V_{j}^n \mathcal{Q}\left(V_{m}^n, V_{n}^n\right)\right\rangle_{\mathbf{v}} K_{m}^n K_{n}^n, \quad j =1,\dots,r. 
\eal
\end{equation}
\item Compute the QR decomposition of $K_j^{n+1}=\sum_{i=1}^rX_i^{n+1}S_{ij}^{(1)}$ to obtain updated $X_i^{n+1}$ and $S_{ij}^{(1)}$. 
\end{enumerate}
The overall arithmetic complexity of this step is $\mathcal{O}\left( r^3N_{\vv}^d+r^3N_{\vx}^d+r^2MN_{\vv}^d\log N_{\vv}\right)$ (suppose the algorithm in \cite{mouhot2006fast} is used for evaluating the collision operator).
\item {\bf $S$-step}.
\begin{enumerate}
\item Perform the forward Euler step in (\ref{eqn:S-step detailed}) to obtain $S_{ij}^{(2)}$:
\begin{equation}
\bal
S_{ij}^{(2)}&=S_{ij}^{(1)} +\Delta t\sum\limits_{l=1}^r\langle \vv V_j^n V_l^n\rangle_{\vv}\cdot \langle X_i^{n+1} \nabla_{\vx}K_l^{n+1}\rangle_{\vx}\\
&-\Delta t\sum\limits_{m,n=1}^r\langle V_j^n \mathcal{Q}(V_m^n,V_n^n)\rangle_{\vv}\suml_{l=1}^{r}\left(\suml_{k=1}^r\left(\langle X_i^{n+1} X_k^{n+1}X_l^{n+1}\rangle_{\vx}S_{km}^{(1)}\right)S_{ln}^{(1)}\right), \quad i,j=1,\dots,r.
\eal
\end{equation}
Since some of the quantities have been computed in the $K$-step, they can be reused in this step, for example, the term $\mathcal{Q}(V_m^n,V_n^n)$. Note that we changed the second term on the right hand side such that it uses $\nabla_{\vx}K_j^{n+1}$ rather than $\nabla_{\vx}X_j^{n+1}$. This is crucial because we have only available the boundary condition expressed in terms of $K_j^{n+1}$ as seen in Section \ref{sec:boundary}.
\end{enumerate}
The overall arithmetic complexity of this step is $\mathcal{O}\left(r^3N_{\vx}^d+r^4\right)$.
\item {\bf $L$-step}.
\begin{enumerate}
\item Construct $L_i^n= \sum_{j=1}^r S_{ij}^{(2)}V_j^n$ and  $\tilde{K}^{n+1}_j = \sum_{i = 1}^r X_i^{n+1}S_{ij}^{(2)}$.
\item Perform the forward Euler step in (\ref{eqn:L-step detailed}) to obtain $L_{i}^{n+1}$:
\begin{equation}
\bal
L_i^{n+1} &=L_i^n-\Delta t\sum\limits_{l=1}^r \vv\cdot \langle X_i^{n+1} \nabla_{\vx} \tilde{K}_l^{n+1}\rangle_{\vx} V_l^n\\
&+\Delta t\sum\limits_{p,q = 1}^r\mathcal{Q}(V_p^n,V_q^n)\suml_{n = 1}^{r}\left(\suml_{m = 1}^r\left(\langle X_i^{n+1} X_m^{n+1}X_n^{n+1}\rangle_{\vx} S_{mp}^{(2)}\right)S_{nq}^{(2)}\right), \quad i=1,\dots,r.
\eal
\end{equation}
The term involving the collision operator is rearranged so that the previously computed values can be reused. For the same reason as in the $S$-step, $\nabla_{\vx} \tilde{K}_j^{n+1}$ is introduced to avoid using $\nabla_{\vx}X_l^{n+1}$.
\item Compute the QR decomposition of $L_i^{n+1}= \sum_{j= 1}^r S_{ij}^{n+1}V_j^{n+1}$ to obtain updated $V_j^{n+1}$ and $S_{ij}^{n+1}$.
\end{enumerate}
The overall arithmetic complexity of this step is $\mathcal{O}\left(r^2N_{\vx}^d+r^3N_{\vv}^d+r^4\right)$.
\end{enumerate}
To simplify the notation, we treat $\vx$, $\vv$ as the continuous variables in the above presentation. The discretization in $\vx$ and $\vv$ can be added straightforwardly following the discussion in Section \ref{velocity discretization} and Section \ref{space discretization}. The inner products $\langle \ \cdot \ \rangle_{\vv}$, $\langle \ \cdot \ \rangle_{\vx}$ are evaluated using the midpoint rule at the discrete velocity and spatial grid points.

If $r$ is small, the computational complexity of the above algorithm will be dominated by the evaluation of the collision operator $\mathcal{O}(r^2MN_{\vv}^d\log N_{\vv})$, which can be much more efficient than the full tensor method whose complexity is $\mathcal{O}(N_{\vx}^dMN_{\vv}^{d}\log N_{\vv})$.

\section{An adaptive dynamical low rank method}
\label{sec:adaptive}

The dynamical low rank method introduced in the last section uses a fixed rank $r$ throughout the entire time evolution. This turns out to be a bad strategy when solving the stationary Boltzmann equation subject to inflow or Maxwell diffusive boundary conditions. The reason is two-fold: 1) The boundary keeps sending new information to the interior of the domain so that the basis $X_i$, $S_{ij}$, $V_j$ initialized according to the initial condition is not sufficient to capture the solution at later time. Thus new basis needs to be injected to the solution over time. 2) For many benchmark tests of the Boltzmann equation, the steady state solutions are often low rank (see Section \ref{sec:lowrank} for a partial justification). Therefore, keeping adding basis without dropping anything would unnecessarily increase the computational cost. In this section, we provide an adaptive strategy to add and delete basis during the time evolution of a dynamical low rank method.

\subsection{Adding basis from the boundary}

Assume that at the boundary $\vx \in \partial \Omega_{\vx}$, $f(\vx,\vv,t)$ is given by 
\begin{equation}
f(\vx,\vv,t)=f^b(\vx,\vv,t), 
\end{equation}
where $f^b(\vx,\vv,t)$ is defined in (\ref{eqn:full bdy general}).

 Since the function $f^b(\vx,\vv,t)$ does not necessarily belong to the space spanned by $\{V_j\}_{1 \le j \le r}$, using a fixed set of basis will result in information loss.

We can fix this problem by explicitly adding boundary conditions as basis into $\{V_j\}_{1\le j \le r}$. For example, at the beginning of time step $t^n$, suppose in the fully discrete scheme there are $N_{bx}$ spatial points on the boundary $\p \Omega_{\vx}$, $N_{\vv}^d$ velocity points over the velocity space $\Omega_{\vv}$ and $N_{\vx}^d$ spatial points over the physical space $\Omega_{\vx}$. We can represent the function $f^b(\vx,\vv,t^n)$ using a matrix $F_b \in \mathbb{R}^{N_{bx} \times N_{\vv}^d}$. We then proceed as follows:

\begin{enumerate}
	\item Compute the SVD of $F_b$ to obtain $F_b = U_b\Sigma_bQ_b^T$ where $U_b, Q_b$ are orthonormal and $\Sigma_b$ is diagonal with descending singular values.
	\item Drop singular values in $\Sigma_b$ that are smaller than $10^{-10}$. Suppose there are $r_b$ singular values remaining, set $\bar{Q}_b = Q_b(:,1:r_b) \in \mathbb{R}^{N_{\vv}^d \times r_b}$.
	\item Concatenate a random matrix $X_h \in \mathbb{R}^{N_{\vx}^d \times r_b}$ to $X^n$, $\bar{Q}_b$ to $V^n$ and extend $S^n$ with zero padding:
	\be
	\widehat{X} = [X^n,X_h] \in \mathbb{R}^{N_{\vx}^d \times (r+r_b)}, \quad \widehat{S} = \bmat{S^n& \\&\textbf{0}} \in \mathbb{R}^{(r+r_b) \times (r+r_b)}, \quad \widehat{V} = [V^n,\bar{Q}_b] \in \mathbb{R}^{N_{\vv}^d \times (r+r_b)}.
	\ee
	\item Perform the QR decomposition of $\widehat{X}$ and $\widehat{V}$ to orthonormalize new basis as $\widehat{X} = X_qS_x$ and $\widehat{V} = V_qS_v$. Set $S_q = S_x\widehat{S}S_v^T$.
\end{enumerate}
Then $(X_q,S_q,V_q)$ are the new basis and we proceed as in Section \ref{scheme:first order}. If $f^b(\vx,\vv,t) = f^b(\vv,t)$ is spatially homogeneous, we can directly start at step 3 and concatenate $F_b$ to $V^n$.


\subsection{Dropping basis adaptively}

To avoid the rank accumulation from the above procedure, we can decrease the rank $r$ by dropping some small singular values of matrix $(S_{ij})_{1\le i,j \le r}$. 

At the end of time step $t^n$ as described in Section \ref{scheme:first order}, we proceed as follows to adjust the rank:
\begin{enumerate}
	\item Compute the SVD of $S^{n+1} = (S_{ij}^{n+1})_{1 \le i,j \le r}$ to obtain $S^{n+1} = U\Sigma Q^T$, where $U,Q \in \mathbb{R}^{r\times r}$ are orthonormal and $\Sigma \in \mathbb{R}^{r \times r}$ is diagonal with descending singular values.
	\item Drop singular values in $\Sigma$ that are less than some tolerance ${\bf drop\_tol}$. Suppose there are $r'$ singular values remaining, we set $\bar{U}=U(:,1:r')$, $\bar{\Sigma} =\Sigma(1:r',1:r')$ and  $\bar{Q}= Q(:,1:r')$. Define $\bar{S}^{n+1} = \bar{\Sigma}$.
	\item Update the basis as $[\bar{X}_{1}^{n+1},\bar{X}_{2}^{n+1},\ldots,\bar{X}_{r'}^{n+1}] =[X_{1}^{n+1},X_{2}^{n+1},\ldots,X_{r}^{n+1}]\bar{U} $ and $[\bar{V}_{1}^{n+1},\bar{V}_{2}^{n+1},\ldots,\bar{V}_{r'}^{n+1}] =[V_{1}^{n+1},V_{2}^{n+1},\ldots,V_{r}^{n+1}]\bar{Q} $ where $\{\bar{X}_{i}^{n+1}\}_{i = 1,\ldots,r'}$ and $\{\bar{V}_{i}^{n+1}\}_{i = 1,\ldots,r'}$ are the updated spatial and velocity basis functions respectively.
	\end{enumerate}

${\bf drop\_tol}$ plays an important role in overall computational efficiency and accuracy. Large ${\bf drop\_tol}$ causes low accuracy for some high-rank solutions and small ${\bf drop\_tol}$ suffers from heavy computation by large computational rank. We dynamically choose ${\bf drop\_tol}$ according to the accuracy of the current solution. More details are given in Section~\ref{subsec:convergence criterion}.

\section{Normal shock problem and low rank property of the solution}
\label{sec:lowrank}

Generally speaking, it is hard to predict or analyze the rank of the solution to the Boltzmann equation due to its highly nonlinear structure. As such, the dynamical low rank method introduced above is really like a black box solver since one cannot tell in advance the rank of the solution until the actual simulation is run. If the rank turns out to be high, the method becomes slow and might not be competitive to the full tensor method. Nevertheless, in this section we identify a class of problems whose solutions are indeed low rank so that we have confidence about the efficiency of the low rank method.

The normal shock problem \cite{cercignani2000rarefied} is a classical benchmark test in rarefied gas dynamics and has been used to validate all kinds of numerical methods for the nonlinear Boltzmann equation. Consider a plane shock wave perpendicular to a flow. The flow is in the $x_1$ direction. The gas is uniform at upstream infinity ($x_1 \rightarrow - \infty$) and downstream infinity ($x_1 \rightarrow + \infty$) and the whole flow is stationary. We are interested in the shock profile developed in this setup with various Mach numbers.

The governing equation is the following 1D stationary Boltzmann equation:
\be
v_1 \p_{x_1} f = \mathcal{Q}(f,f),
\label{eqn:normalshock}
\ee
with boundary condition
\be
\bal
&\lim\limits_{x_1 \rightarrow -\infty} f(x_1,\vv) = f_L(\vv) = \mathcal{M}(\rho_L,\vu_L,T_L)(\vv) = \frac{\rho_L}{(2\pi R T_L)^{d/2}}\exp\left(-\frac{(v_1-u_L)^2+v_2^2+...+v_d^2}{2RT_L}\right),\\
&\lim\limits_{x_1 \rightarrow + \infty} f(x_1,\vv) = f_R(\vv) = \mathcal{M}(\rho_R,\vu_R,T_R)(\vv)=\frac{\rho_R}{(2\pi R T_R)^{d/2}}\exp\left(-\frac{(v_1-u_R)^2+v_2^2+...+v_d^2}{2RT_R}\right),
\label{eqn:normalshockboundary}
\eal
\ee
where $\mathcal{M}(\rho,\vu,T)$ is the Maxwellian distribution; $(\rho_L,\vu_L,T_L)$ and $(\rho_R,\vu_R,T_R)$ are the density, bulk velocity and temperature of the upstream and downstream flows; and $R$ is the gas constant.

The net flow of mass, momentum and energy into the shock must be equal to the ones out of the shock:
\be
\int v_1f_L(\vv)  \bmat{1\\v_1 \\  \vv^2}\, {\rm d} \vv = \int v_1f_R(\vv) \bmat{1\\v_1 \\ \vv^2}\,{\rm d} \vv.
\label{eqn:normal_shock_conservation1}
\ee
Rewriting \cref{eqn:normal_shock_conservation1} in terms of macroscopic quantities $\rho_{L,R}$, $u_{L,R}$ and $T_{L,R}$, we have the following Rankine-Hugoniot relations
\be
\bal
&\rho_Lu_L = \rho_Ru_R,\\
&\rho_Lu_L^2 + \rho_LRT_L = \rho_Ru_R^2 + \rho_RRT_R, \\
&\rho_Lu_L\left( u_L^2+ (d+2)R T_L\right) = \rho_Ru_R\left(  u_R^2+ (d+2) R T_R\right).  
\label{eqn:normal_shock_conservation2}
\eal
\ee
Given the upstream quantities $(\rho_L, u_L, T_L)$ and using the upstream flow Mach number defined by 
\be
M_L = \frac{u_L}{(\gamma RT_L)^{\fl 1 2}}, \quad \gamma = \fl {d+2} {d},
\ee
we can solve \cref{eqn:normal_shock_conservation2} to obtain
\be \label{relation}
\rho_R = \rho_L \frac{(d+1)M_L^2}{M_L^2+d}, \quad u_R = u_L \frac{M_L^2+d}{(d+1)M_L^2}, \quad T_R = T_L \frac{((d+2)M_L^2-1)(M_L^2+d)}{(d+1)^2M_L^2}.
\ee

In the following, we consider two scenarios where one can obtain some low rank approximation to the solutions of (\ref{eqn:normalshock})-(\ref{eqn:normalshockboundary}).

\subsection{Weak shock wave: \texorpdfstring{$M_L = \mathcal{O}(1)$}{Lg}}

When $M_L=1$, it is clear from (\ref{relation}) that there will be no jump hence no shock. When $M_L=\mathcal{O}(1)$ but bigger than $1$, a weak shock will be developed. We assume 
\begin{equation}
M_L = 1+\varepsilon, 
\end{equation}
where $\varepsilon$ is a small parameter. In fact, $\varepsilon$ is on the same order of the mean free path \cite{Ohwada93}. We then rescale $x_1$ according to $\tilde{x}_1 = \varepsilon x_{1}$. The \cref{eqn:normalshock} thus becomes
\be
v_1 \p_{\tilde{x}_1} f = \frac{1}{\varepsilon} \mathcal{Q}(f,f).
\label{eqn:rescaled normalshock}
\ee
On the other hand, we can see from (\ref{relation}) that the macroscopic quantities of upstream flow and downstream flow are very close:
\be
\bal
&\frac{\rho_R}{\rho_L}  = 1+\frac{d(M_L^2-1)}{M_L^2+d}=1+\mathcal{O}(\varepsilon),\\
&\frac{u_R}{u_L} = 1-\frac{d(M_L^2-1)}{(d+1)M_L^2}=1+\mathcal{O}(\varepsilon),\\
&\frac{T_R}{T_L} =1+\frac{(d+1)(M_L^4-1)+(M_L^2-1)^2}{(d+1)^2M_L^2}= 1+\mathcal{O}(\varepsilon).
\eal
\ee
Hence
\be
\frac{f_R}{f_L} = 1+\mathcal{O}(\varepsilon).
\ee
Therefore, it is reasonable to assume
\be
f(\tilde{x}_1,\vv) = f_L(\vv)+\varepsilon f_1(\tilde{x}_1,\vv)+ \mathcal{O}(\varepsilon^2),
\label{eqn:hibert expansion form}
\ee
where $f_1(\tilde{x}_1,\vv)$ is yet to be determined.

The rest of the analysis is similar to the Hilbert expansion. Substituting (\ref{eqn:hibert expansion form}) into (\ref{eqn:rescaled normalshock}) and matching orders, we obtain at order $\mathcal{O}(\varepsilon)$:
\be
\mathcal{Q}(f_1,f_L) + \mathcal{Q}(f_L,f_1) = v_1\p_{\tilde{x}_1}f_L(\vv)\equiv 0.
\label{eqn:hilbert expansion}
\ee
Using the linearized Boltzmann collision operator \cite{Cercignani} defined by
\be
L_{\mathcal{M}}(f) := \frac 1 {\mathcal{M}} \left(\mathcal{Q}(\mathcal{M},\mathcal{M}f)+\mathcal{Q}(\mathcal{M}f,\mathcal{M})\right), \quad \mathcal{M} \text{ is a Maxwellian},
\ee
we can write (\ref{eqn:hilbert expansion}) as
\be
L_{f_L}\left(\frac{f_1}{f_L}\right)(\tilde{x}_1,\vv) = 0.
\label{eqn:f0hilbert}
\ee
The kernel property of $L_{\mathcal{M}}$ implies that $\frac{f_1}{f_L}$ must be a linear combination of collision invariants 1, $\vv$, $|\vv|^2$ and we may write
\be
{f_1}(\tilde{x}_1,\vv) = f_L(\vv)\left(a(\tilde{x}_1)+ {\bf b}(\tilde{x}_1)\cdot \vv+ c(\tilde{x}_1)|\vv|^2\right),
\ee
where $a$, $\bf b$ and $c$ are functions of $x_1$ only. Together with (\ref{eqn:hibert expansion form}), we have
\be
f(\tilde{x}_1,\vv) = f_L(\vv)(1+\varepsilon a(\tilde{x}_1)+\varepsilon {\bf b}(\tilde{x}_1)\cdot \vv+ \varepsilon c(\tilde{x}_1)|\vv|^2 )+ \mathcal{O}(\varepsilon^2).
\ee
Therefore, up to order $O(\varepsilon)$, the solution $f(\tilde{x}_1,\vv)$ is a low rank separated function in $\tilde{x}_1$ and $\vv$.

We mention that the derivation of $\mathcal{O}(\varepsilon)$ term does not require specific properties of the collision kernel $B$. One can continue this process to derive $\mathcal{O}(\varepsilon^2)$ term, which is a low rank function as well and depends on the kernel $B$, see \cite{Ohwada93} for details.

\subsection{Strong shock wave: \texorpdfstring{$M_L \rightarrow \infty$}{Lg}}

When $M_L$ is very large, a strong shock wave will develop and one cannot hope for any asymptotic expansion as in the previous subsection. Over the years, people have tried to find various approximations to the solution in this regime and it turns out many heuristic solutions match well with the experiments, yet are low rank \cite{cercignani2000rarefied,Harris}. Here we present one such approximation due to Mott-Smith, who obtained the first solution of Boltzmann's equation for the shock structure problem in 1951. More sophisticated approximations exist but they more or less follow a similar idea as Mott-Smith.

The starting point is a bimodal distribution (and low rank) approximation of $f$ as
\be
f(x_1,\vv) = a_1(x_1)f_L(\vv) + a_2(x_1)f_R(\vv).
\label{approximation: strong shock wave}
\ee
To satisfy the Rankine-Hugoniot equations, we must have $a_1(x_1)+a_2(x_1) \equiv 1$. We thus write $a(x_1) = a_1(x_1)$ and $a_2(x_1) = 1-a(x_1)$. In order to determine $a(x_1)$, one additional condition is needed. The simplest way is to enforce the moment equation by multiplying \cref{eqn:normalshock} by $\int \cdot \ v_1^2 \rd{\vv}$:
\be
\int v_1^3 \p_{x_1} f \, {\rm d} \vv=  \int v_1^2\mathcal{Q}(f,f) \, {\rm d} \vv,
\ee
which reduces to 
\be \label{aa}
a'(x_1)\left(\rho_Lu_L(u_L^2+3RT_L)-\rho_Ru_R(u_R^2+3RT_R)\right) = {\alpha} a(x_1)(1-a(x_1)),
\ee
with
\be
\bal
\alpha &= \int v_1^2 \left(\mathcal{Q}(f_L,f_R)+\mathcal{Q}(f_R,f_L)\right) \, {\rm d} \vv.
\label{eqn:alpha}
\eal
\ee
Using (\ref{eqn:normal_shock_conservation2}), (\ref{aa}) can be further simplified to
\be
(d-1)\rho_Lu_L R(T_L-T_R)a'(x_1)  = -{\alpha} a(x_1)(1-a(x_1)).
\ee
This equation easily integrates to 
\be
a(x_1) = \frac {1} {\exp(\beta x_1)+1}, \quad \beta = \frac{\alpha}{(d-1)\rho_Lu_LR(T_L-T_R)}.
\ee

Therefore, we have found a closed form solution in the form of (\ref{approximation: strong shock wave}). Note that to evaluate $\alpha$, we need to make use of specific properties of the collision kernel $B$. Accordingly, we can see that the spatial change in $\rho$ across the shock wave increases with increasing Mach number $M_L$ of the upstream:
\be
\frac{\rho(x_1)}{\rho_L} = \frac{1+ \frac{(d+1)M_L^2}{M_L^2+d}\exp(\beta x_1)}{1+\exp(\beta x_1)}.
\ee


\section{Numerical examples}
\label{sec:num}

In this section, we evaluate the accuracy and efficiency of the proposed dynamical low rank method by several classical benchmark tests in rarefied gas dynamics, including normal shock wave (1D), Fourier flow (1D), lid driven cavity flow (2D), and thermally driven cavity flow (2D). All these examples concern the steady-state solution of the Boltzmann equation and we use the first order method in both time and space as described in Section \ref{sec:dynamical}, and Fourier spectral method for 2D Maxwell molecules \cite{mouhot2006fast} to evaluate the collision operator. The results are compared with full tensor method using the same discretization.

\subsection{Convergence criterion}
\label{subsec:convergence criterion}
Unlike time dependent problems, we need a proper stopping criterion for solving the steady state solutions.

For the full tensor method, we define the error as
\begin{equation}
\begin{aligned}
\text{err}_\text{\text{full tensor}} &= \lVert f_{\text{full tensor}}^{n+1}-f_{\text{full tensor}}^{n} \rVert_{L^2} = \left\langle f_{\text{full tensor}}^{n+1}-f_{\text{full tensor}}^{n},f_{\text{full tensor}}^{n+1}-f_{\text{full tensor}}^{n}\right\rangle_{\vx,\vv}^{\fl 1 2}.
\label{err full grid}
\end{aligned}
\end{equation}

For the low rank method, we define the error similarly as
\begin{equation}
\text{err}_{\text{low}\ \text{rank}} = \lVert f_{\text{low}\ \text{rank}}^{n+1}-f_{\text{low}\ \text{rank}}^{n} \rVert_{L^2} = \left\langle f_{\text{low}\ \text{rank}}^{n+1}-f_{\text{low}\ \text{rank}}^{n},f_{\text{low}\ \text{rank}}^{n+1}-f_{\text{low}\ \text{rank}}^{n}\right\rangle_{\vx,\vv}^{\fl 1 2},
\label{error:low rank nonada}
\end{equation}
where $f^n_{\text{low rank}}=\sum_{i,j=1}^r X_i^nS_{ij}^nV_j^n$. Rather than reconstructing $f^n_{\text{low rank}}$, the above error term can be broke into three pieces:
\begin{equation}
\begin{aligned}
f_{\text{low}\ \text{rank}}^{n+1}-f_{\text{low}\ \text{rank}}^{n}& =\suml_{i,j=1}^r X_i^{n+1}S_{ij}^{n+1}V_j^{n+1}- \suml_{i,j=1}^rX_i^{n}S_{ij}^{n}V_j^{n}\\
& =  \suml_{j=1}^r \left( K_j^{n+1}-K_j^{n}\right)V_j^{n}+\suml_{i,j=1}^r X_i^{n+1}\left( S_{ij}^{(2)}-S_{ij}^{(1)}\right)V_j^{n}+\suml_{i=1}^rX_i^{n+1}\left(L_i^{n+1}-L_i^{n}\right)\\
&:= \suml_{j=1}^r \dt K_j V_j^{n}+\suml_{i,j=1}^r X_i^{n+1}\dt S_{ij} V_j^{n}+\suml_{i=1}^rX_i^{n+1} \dt L_i
\end{aligned}
\end{equation}
where the notation follows Section \ref{scheme:first order}. By orthogonality of $\{X_i\}_{1 \le i \le r}$ and $\{V_j\}_{1 \le j \le r}$, (\ref{error:low rank nonada}) can be simplified as
\begin{equation}
\begin{aligned}
\text{err}_{\text{low}\ \text{rank}}^2 &= \left\langle f_{\text{low}\ \text{rank}}^{n+1}-f_{\text{low}\ \text{rank}}^n,f_{\text{low}\ \text{rank}}^{n+1}-f_{\text{low}\ \text{rank}}^n \right\rangle_{\vx,\vv}\\
& = \suml_{j=1}^r\left\langle \dt K_j,\dt K_j \right\rangle_{\vx}+\suml_{i,j=1}^r\dt S_{ij}^2+\suml_{i=1}^r \left\langle \dt L_i,\dt L_i\right\rangle_{\vv}+ \text{I}+ \text{II}+\text{III},
\label{err low rank1}
\end{aligned}
\end{equation}
where \text{I}, \text{II} and \text{III} are cross terms:
\begin{equation}
\begin{aligned}
&\text{I} = 2\suml_{i,j= 1}^r  \left\langle \dt K_j,X_i^{n+1}\right\rangle_{\vx}\dt S_{ij}.\\
&\text{II} = 2\suml_{i,j=1}^r \left\langle \dt L_i,V_j^n\right\rangle_{\vv}\dt S_{ij}.\\
&\text{III} = 2\suml_{i,j=1}^r \left\langle \dt K_j,X_i^{n+1}\right\rangle_{\vx}\cdot \left\langle \dt L_i,V_j^n \right \rangle_{\vv}.
\label{err low rank2}
\end{aligned}
\end{equation}
We emphasize that it is crucial to evaluate $\text{err}_{\text{low}\ \text{rank}}$ using (\ref{err low rank1})-(\ref{err low rank2}), since the cost of reconstructing $f^n_{\text{low rank}}$ is $\mathcal{O}(r^2N_{\vx}^dN_{\vv}^d)$ which is comparable to a full tensor method.

In general, we set a fixed convergence tolerance \text{\bf res\_tol} and terminate the time iteration whenever $\text{err}_{\text{low}\ \text{rank}}, \ \text{err}_{\text{full}\ \text{tensor}} \le {\bf \text{\bf res\_tol}}$ for both the full tensor method and low rank method.

For the adaptive low rank method discussed in Section \ref{sec:adaptive}, we have
\begin{equation}
\begin{aligned}
	\abs{\text{err}_{\text{low}\ \text{rank}}-\text{err}_{\text{low}\ \text{rank}}^\text{ada}} \le \|\bar{f}_{\text{low}\ \text{rank}}^{n+1}-f^{n+1}_{\text{low}\ \text{rank}}\|_{L^2} \le (r-r')^{\fl 1 2}\cdot {\bf drop\_tol},
	\label{ada bound}
\end{aligned}
\end{equation}
where $\text{err}_{\text{low}\ \text{rank}}^\text{ada}=\|\bar{f}_{\text{low}\ \text{rank}}^{n+1}-f^n_{\text{low}\ \text{rank}}\|_{L^2}$, $\bar{f}_{\text{low}\ \text{rank}}^{n+1}$ is the solution at the end of time step $t^n$ after adding and removing basis. We dynamically set ${\bf drop\_tol} = c \cdot \text{err}_{\text{low}\ \text{rank}}^\text{ada}$ and control $\text{err}_{\text{low}\ \text{rank}}^\text{ada}$ through
\begin{equation}
\frac{1}{1+c(r-r')^{\fl 1 2}} \text{err}_{\text{low}\ \text{rank}} \le \text{err}_{\text{low}\ \text{rank}}^\text{ada} \le \frac{1}{1-c(r-r')^{\fl 1 2}} \text{err}_{\text{low}\ \text{rank}}.
\end{equation}
In the following tests, we set $c = 0.2$ and always use the adaptive dynamical low rank method with convergence criterion $\text{err}_{\text{low}\ \text{rank}}^\text{ada} \le {\bf \text{\bf res\_tol}}$. 

\subsection{Normal shock wave}

We first consider the normal shock problem (\ref{eqn:normalshock})-(\ref{eqn:normalshockboundary}) with several different Mach numbers $M_L$. We take $R=1$, $d=2$, hence $\gamma=2$, $M_L=\frac{u_L}{(2T_L)^{1/2}}$. In the following, the spatial domain is chosen as $x_1 \in [-30,30]$ with $N_{\vx} = 1000$; and the velocity domain is $(v_1,v_2) \in [-L_{\vv},L_{\vv}]^2$. 

We choose the upstream and downstream condition as
\[
(\rho_L,\ \rho_R) = \left(1,\ \fl{3M_L^2}{M_L^2+2}\right),\quad  (u_L,\ u_R) = \left(\sqrt{2}M_L,\ \fl{\rho_Lu_L}{\rho_R}\right),\quad (T_L,\ T_R) = \left(1,\ \fl{4 M_L^2-1}{3\rho_R}\right),
\]
and the initial condition as
\[
\rho_0(x_1) = \fl{\tanh(\alpha x_1)+1}{2(\rho_R-\rho_L)}+\rho_L,\quad T_0(x_1) = \fl{\tanh(\alpha x_1)+1}{2(T_R-T_L)}+T_L,\quad \vu_0(x_1) = \left(\fl{\tanh(\alpha x_1)+1}{2(u_R-u_L)}+u_L,0\right),
\]
with $\alpha = 0.5$.

When showing the numerical results, we are mainly interested in the macroscopic quantities: density $\rho(x_1)$, bulk velocity $u(x_1)$ (in first dimension) and temperature $T(x_1)$. Their normalized values will be plotted, which are defined by
\[
\widehat{\rho}(x_1) = \frac{\rho(x_1)-\rho_L}{\rho_R-\rho_L}, \quad \widehat{u}(x_1) = \frac{u(x_1)-u_R}{u_L-u_R}, \quad \widehat{T}(x_1) = \frac{T(x_1)-T_L}{T_R-T_L}.
\]

\subsubsection{Weak shock wave: Mach 1.4}

In this subsection we consider Mach number to be $M_L= 1.4$ and set $N_{\vv} = 32$, $L_{\vv} = 13.11$. We set the reference solution $f_{\text{ref}}$ as the solution from the full grid method with convergence criterion $\text{\bf res\_tol} = 4\times 10^{-10}$. 

We check both the full grid method and adaptive low rank method by varying convergence criterion $\text{\bf res\_tol}$. The error is defined as $\lVert f_{\text{ref}}-f_{\text{num}}\rVert_{L^2}$ where $f_{\text{num}}$ is the solution from either the full grid or low rank method. At the same time, we record the computational time needed for both methods to reach the same convergence criterion.

\begin{figure}[!htbp]
\begin{center}
\includegraphics[width=2.8in]{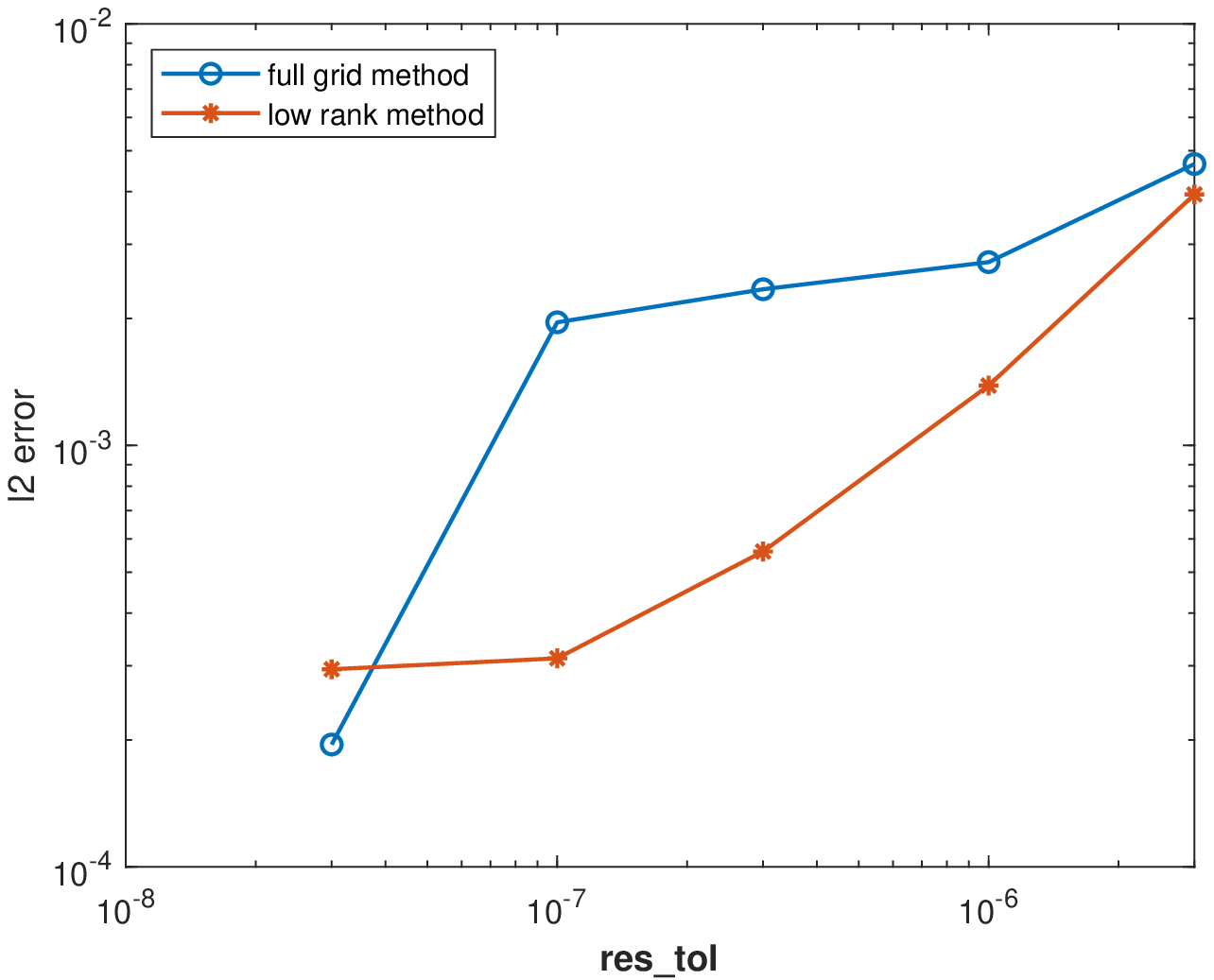}
\includegraphics[width=2.8in]{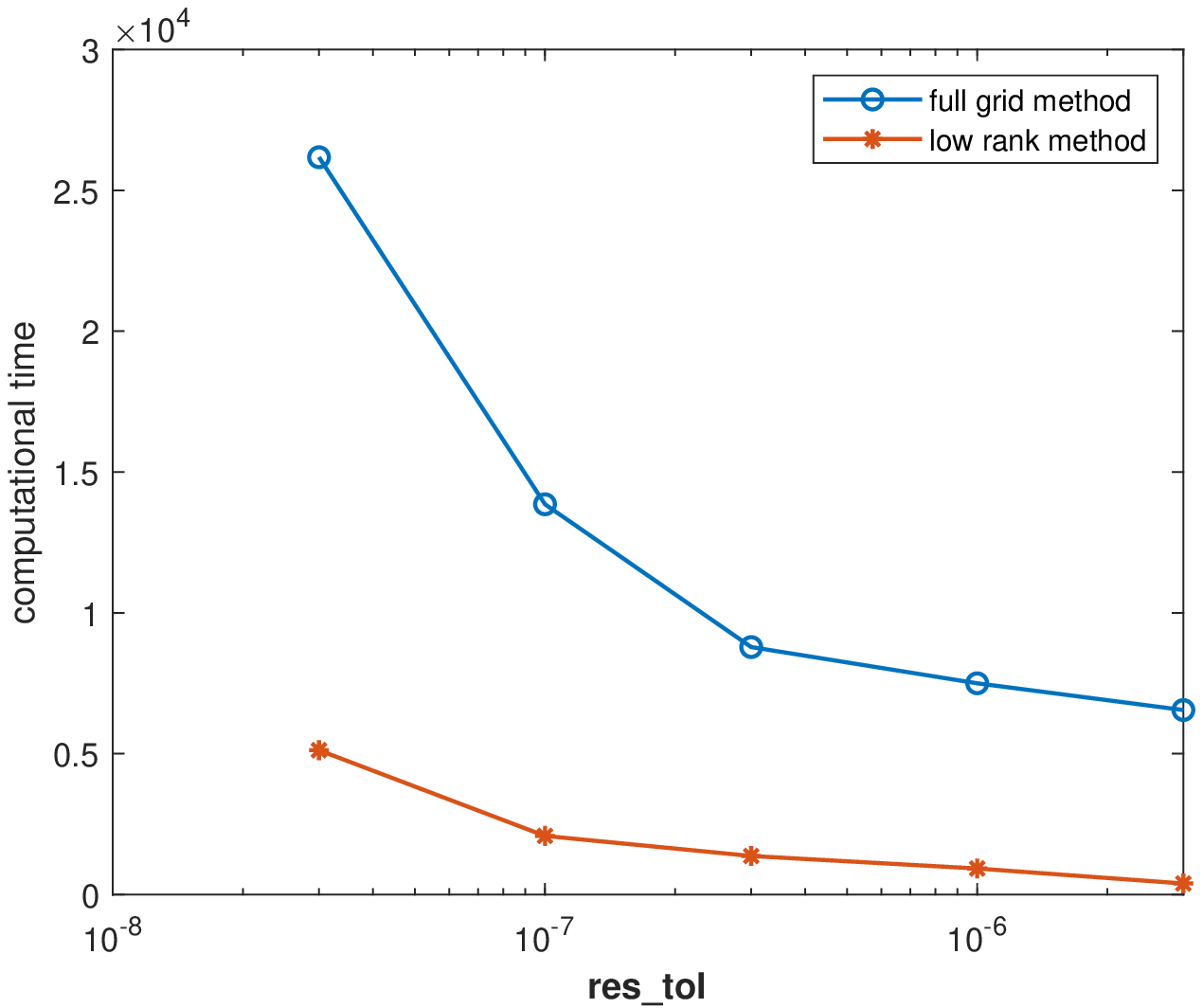}
\caption{Normal shock wave (Mach 1.4). Left: error of the full grid method and the adaptive low rank method for different convergence criterion $\text{\bf res\_tol}$. Right: computational time in seconds for both methods.}
\label{normal shock:mach14 err/time}
\end{center}
\end{figure}

\begin{figure}[!htbp]
\begin{center}
\includegraphics[width=2.8in]{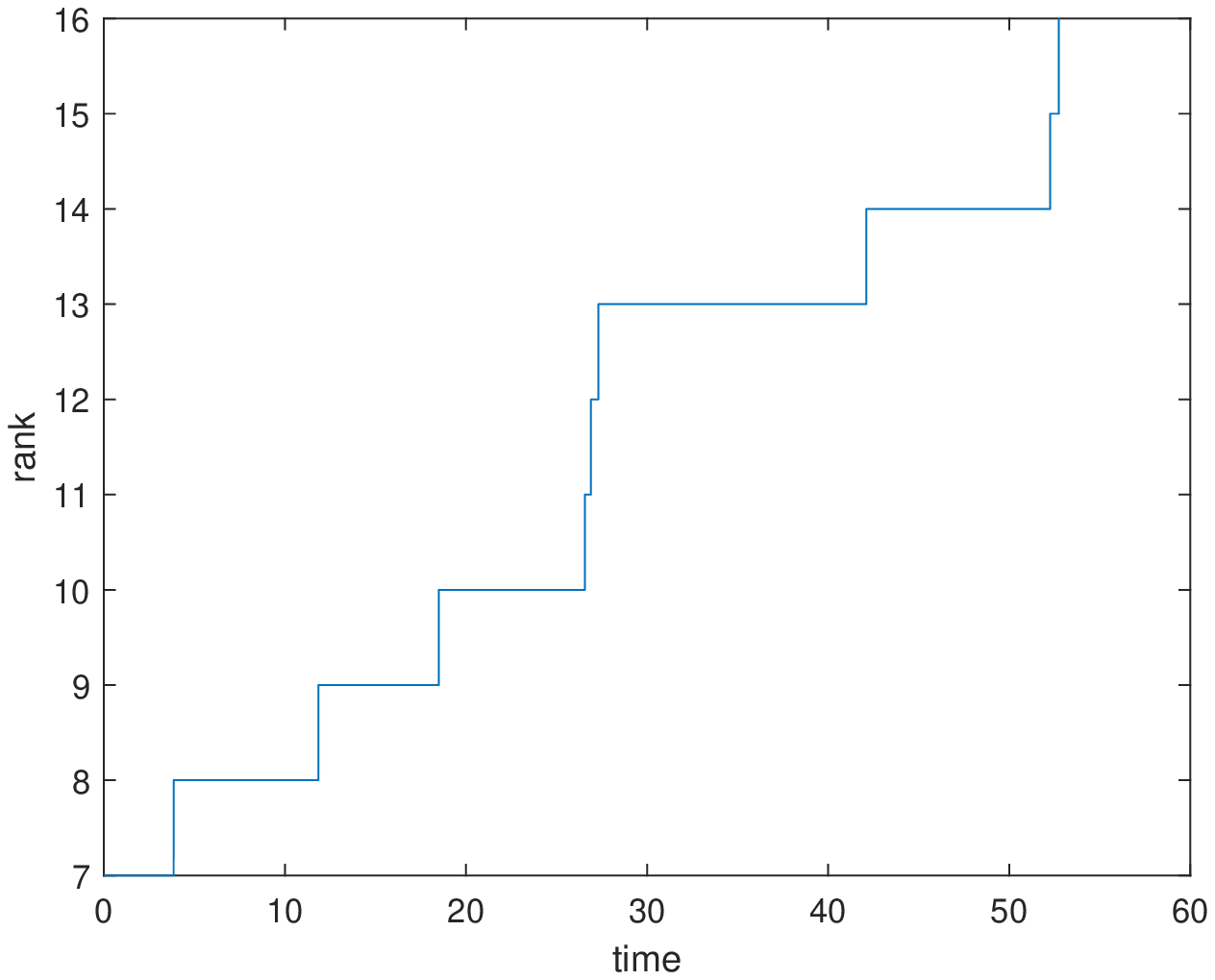}
\includegraphics[width=2.8in]{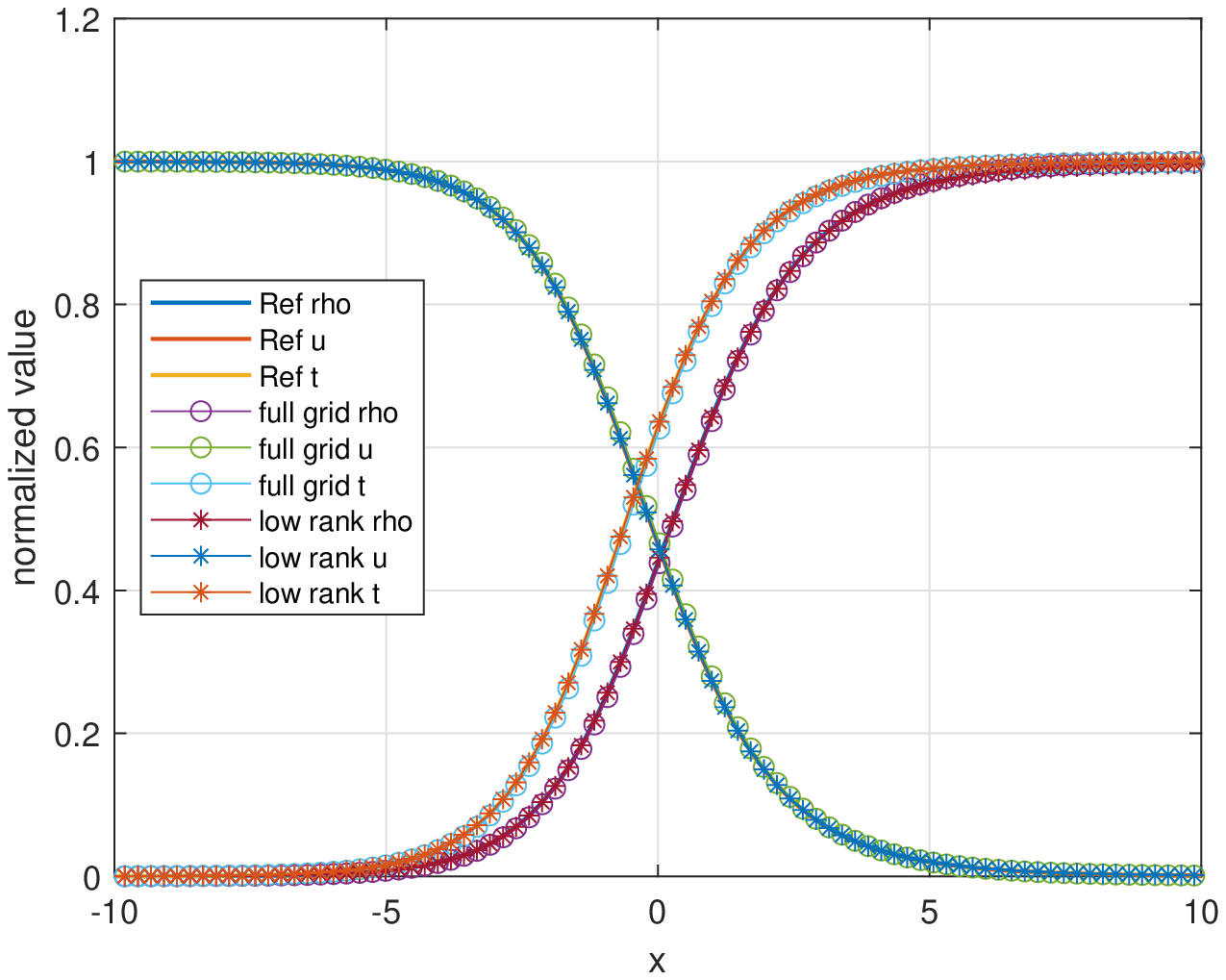}
\caption{Normal shock wave (Mach 1.4). Left: rank evolution of the adaptive low rank method. Right: normalized density, bulk velocity and temperature of the full grid method and the adaptive low rank method using $\text{\bf res\_tol} = 3\times 10^{-7}$, in comparison to the reference solution. }
\label{normal shock:mach14 mac}
\end{center}
\end{figure}

From \cref{normal shock:mach14 err/time}, we can see that the low rank method can achieve the same accuracy much more efficiently compared to the full grid method. From \cref{normal shock:mach14 mac}, we can see that both methods match well with the reference solution. The rank in the adaptive low rank method grows slowly as time evolves and is stabilized to 16 before reaching the convergence criterion.


\subsubsection{Strong shock wave: Mach 3.8 \& Mach 6.5}

In this subsection we consider the strong shock wave with two different Mach numbers $M_L= 3.8$ and $M_L= 6.5$. We compare the full grid method and the adaptive low rank method using the same convergence criterion $\text{\bf res\_tol} = 4.6\times 10^{-7}$.

For the case $M_L = 3.8$, we use $N_{\vv} = 32$ and $L_{\vv}= 20.97$. The full grid method needs 18540 seconds to converge; and the adaptive low rank method needs 7556 seconds to converge. For the case $M_L = 6.5$, we use $N_{\vv} = 48$ and $L_{\vv} = 34.08$. The full grid method needs 44379 seconds to converge; and the adaptive low rank method needs 16157 seconds to converge. 

\begin{figure}[!htbp]
\begin{center}
\includegraphics[width=2.8in]{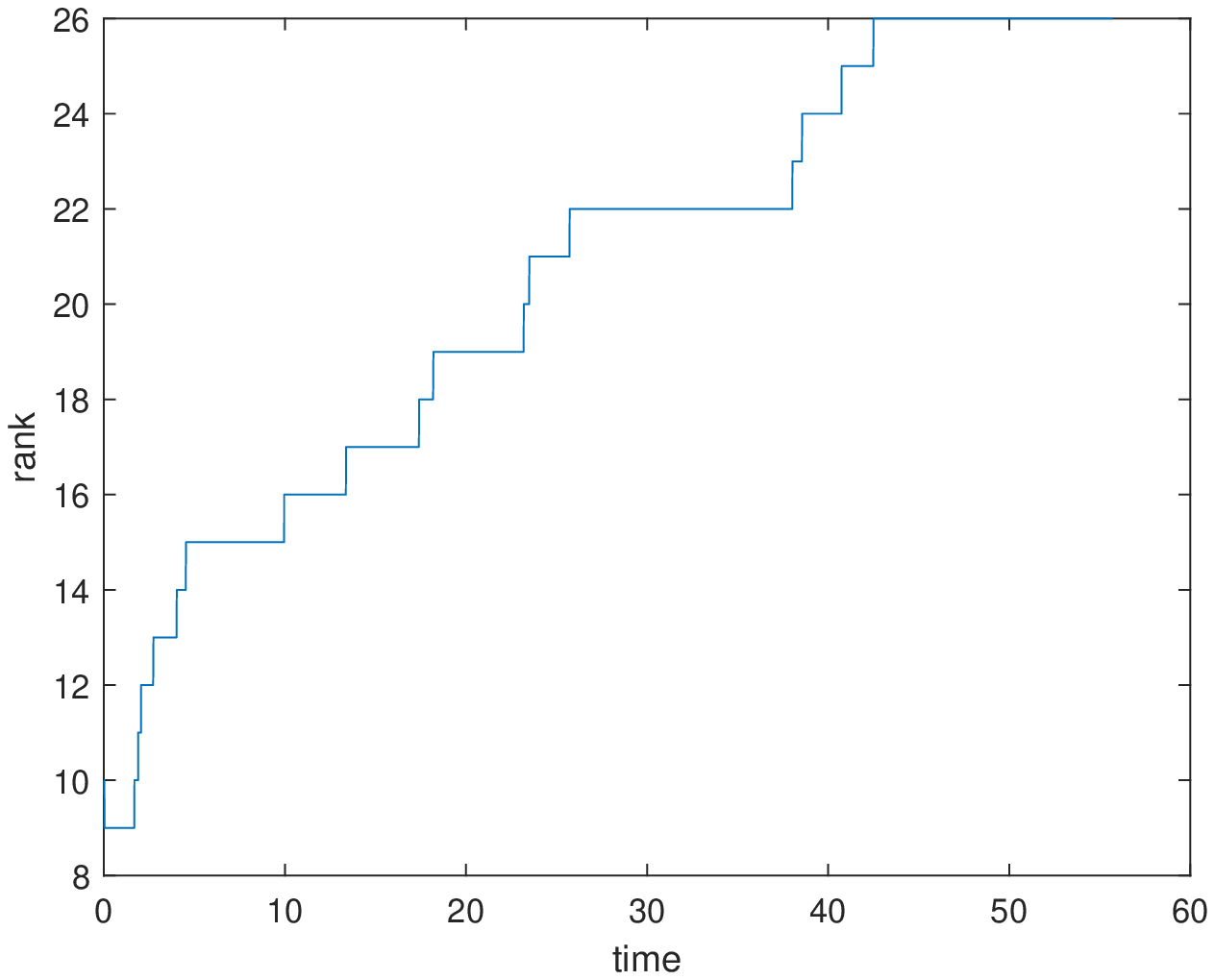}
\includegraphics[width=2.8in]{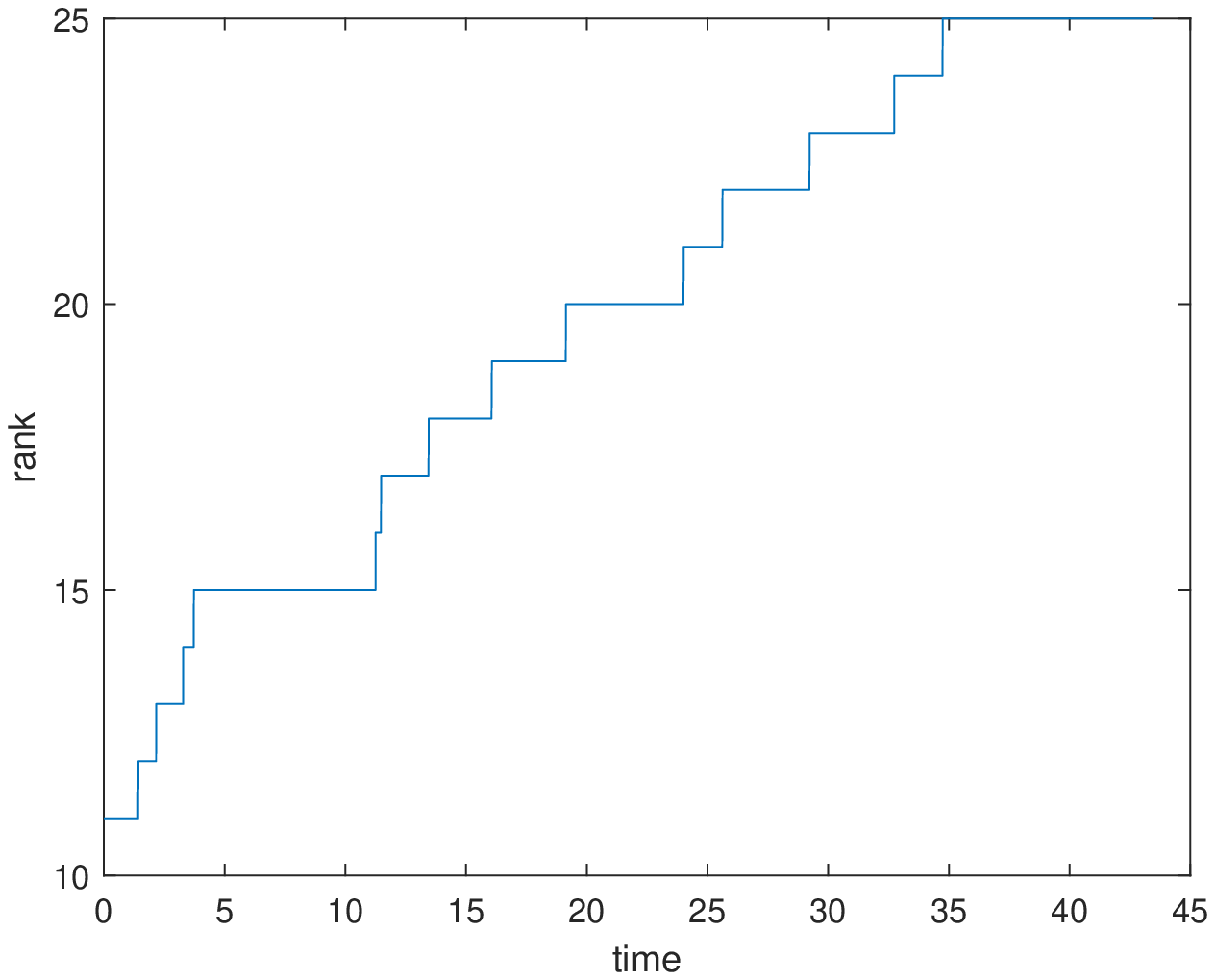}
\includegraphics[width=2.8in]{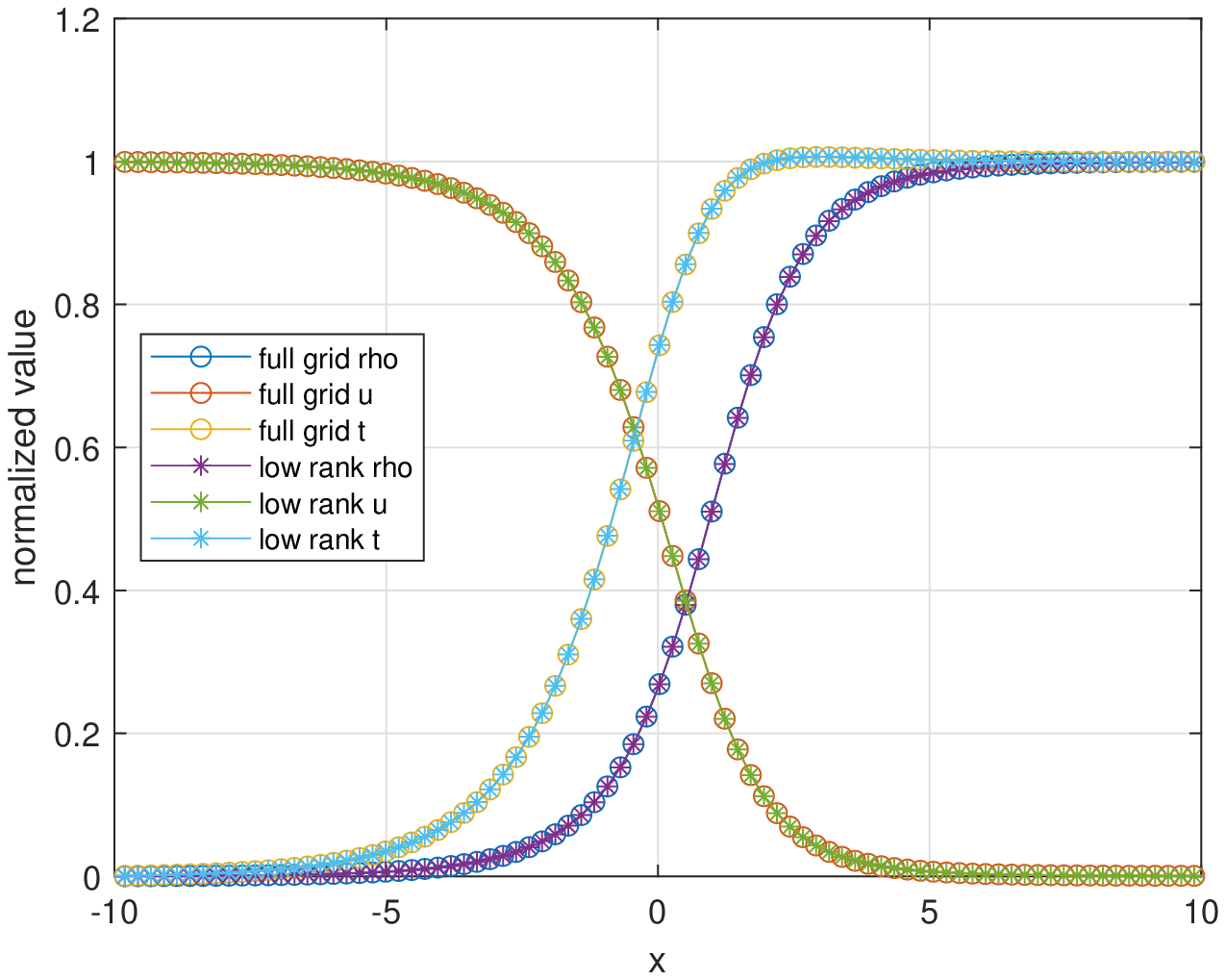}
\includegraphics[width=2.8in]{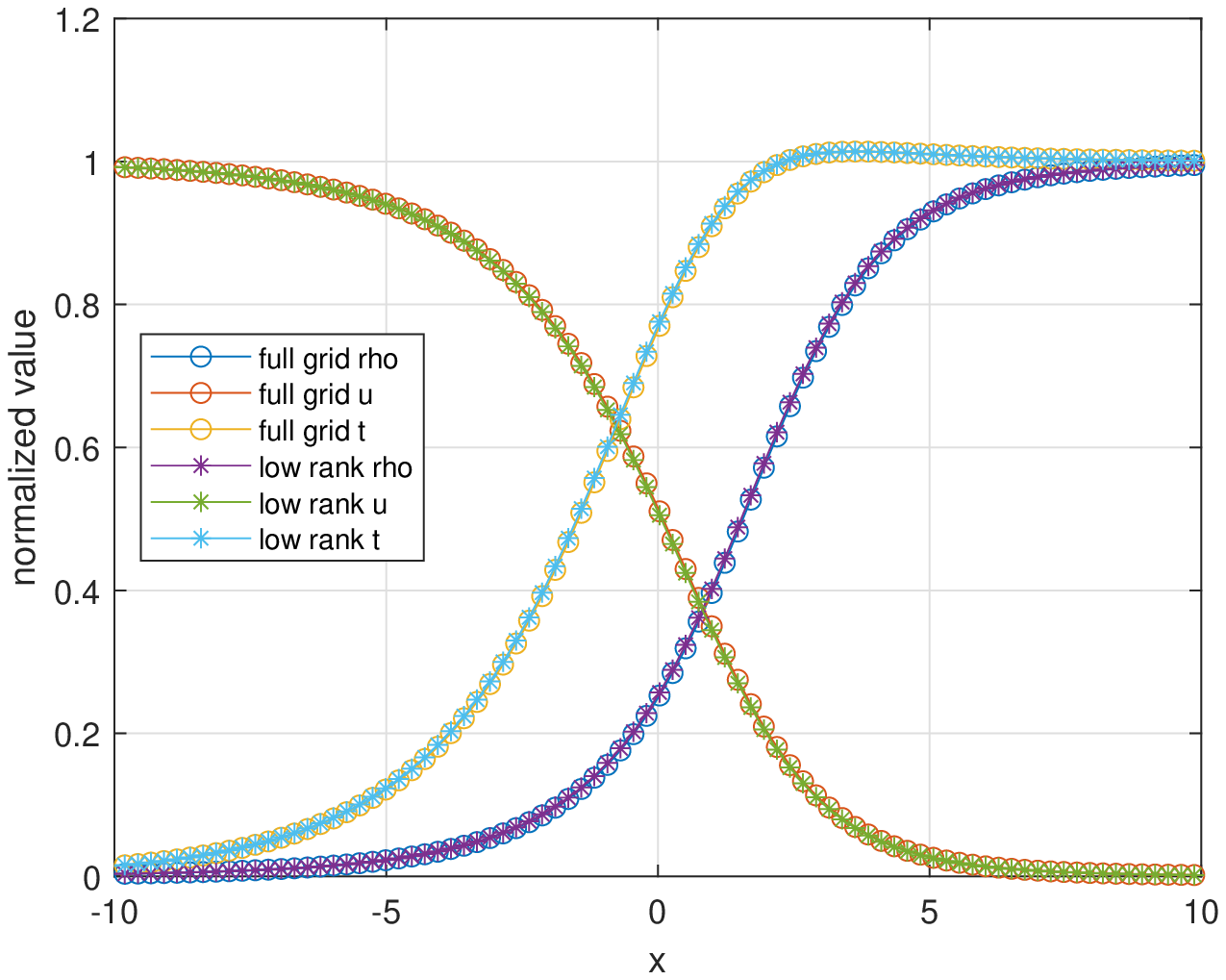}
\caption{Normal shock wave (Mach 3.8 \& Mach 6.5). Top: rank evolution of the adaptive low rank method with Mach 3.8 (Left) and Mach 6.5 (Right); Bottom: normalized density, bulk velocity and temperature of the full grid method and the adaptive low rank method using $\text{\bf res\_tol} = 4.6\times 10^{-7}$ with Mach 3.8 (Left) and Mach 6.5 (Right).}
\label{normal shock: large mach}
\end{center}{}
\end{figure}

The results of both cases are reported in \cref{normal shock: large mach}. The full grid method and the adaptive low rank method match well. On the other hand, the rank in the adaptive low rank method behaves similarly as in the weak shock wave: the numerical rank is a bit higher but still quite low rank and stabilized before reaching the convergence criterion.


\subsection{Fourier flow}

We next consider a Fourier heat transfer problem. The spatial domain is 1D: $x_1 \in [0,2]$ with $N_{\vx} = 200$; and the velocity domain is 2D: $(v_1,v_2) \in [-L_{\vv},L_{\vv}]^2$ with $L_{\vv} = 7.86$ and $N_{\vv} = 32$. The Maxwell diffusive boundary condition is assumed at $x_1 = 0$ with wall quantities $\vu_{w} = (0,0)$, $T_{w} = 1$ and $x_1 = 2$ with $\vu_{w} = (0,0)$, $T_{w} = 1.2$. For the initial condition, we use a spatially homogeneous Maxwellian with $\rho_0 = 1$, $\vu_0 = (0,0)$ and $T_0 = 1$. 

\begin{figure}[!htbp]
\begin{center}
\includegraphics[width=2.8in]{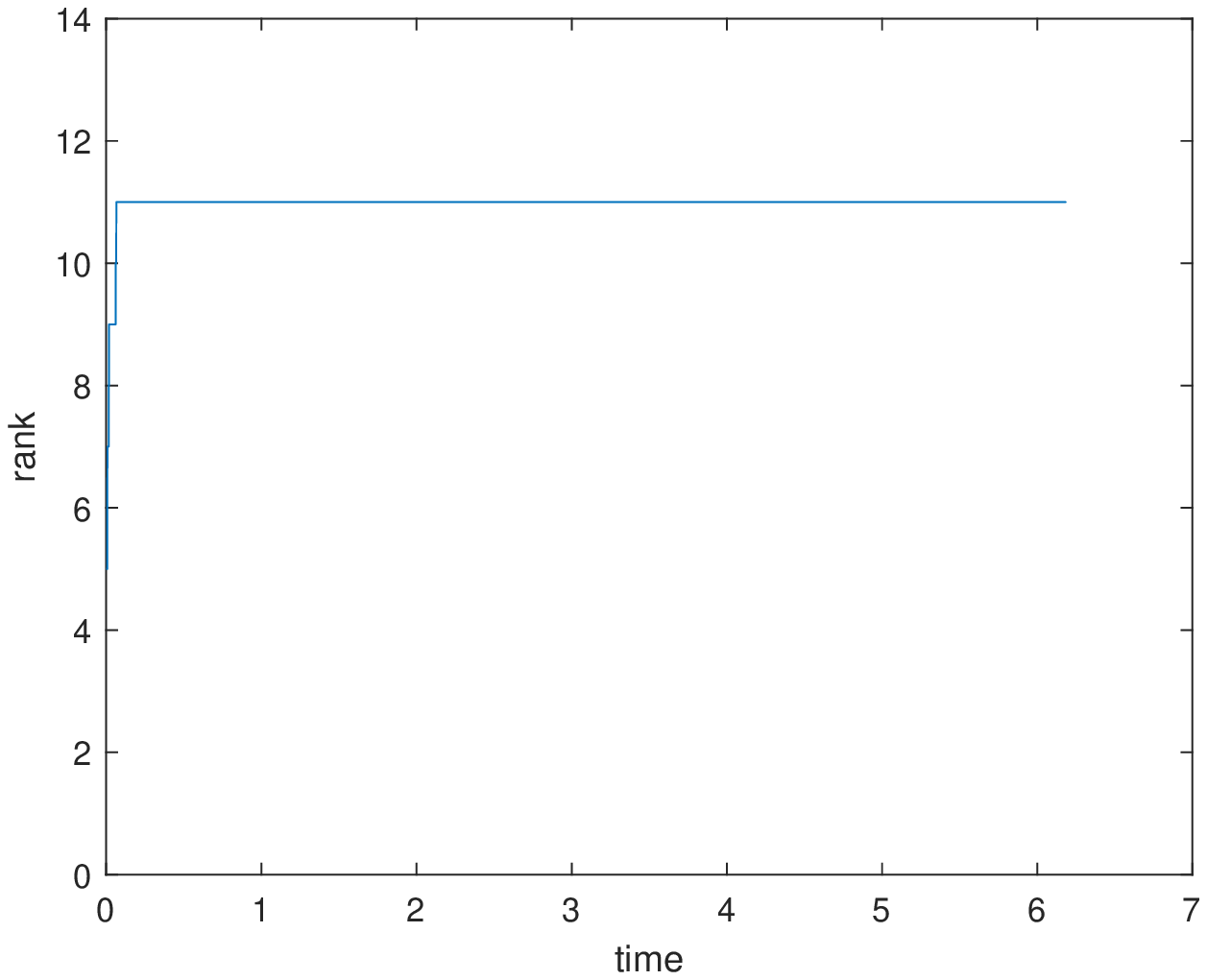}
\includegraphics[width=2.8in]{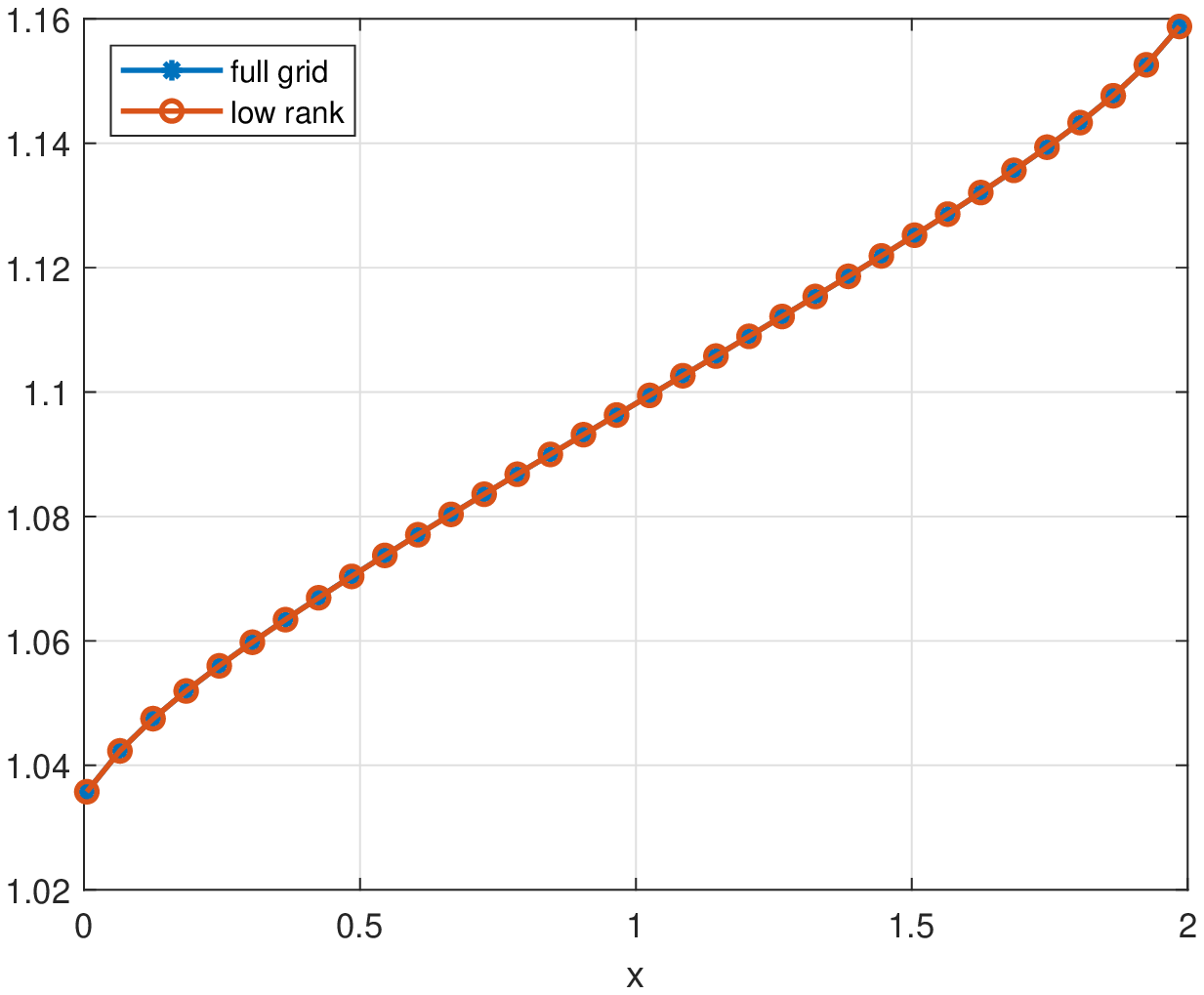}
\caption{Fourier flow. Left: rank evolution in the adaptive low rank method; Right: temperature profile of the full grid method and the adaptive low rank method using $\text{\bf res\_tol} = 2\times 10^{-7}$.}
\label{fourier flow temp}
\end{center}
\end{figure}

The convergence criterion is set as $\text{\bf res\_tol} = 2\times 10^{-7}$ for both the full grid method and the adaptive low rank method. For the full grid method, we need 925 seconds to reach convergence, while for the low rank method, we only need 509 seconds. The temperature profile as shown in \cref{fourier flow temp} matches well for both methods. Furthermore, we can see that the numerical rank in the adaptive low rank method is stabilized to 11 in a very short time.

\subsection{Lid driven cavity flow}

We now consider the 2D lid driven cavity flow problem. The spatial domain is rectangular $(x_1,x_2) \in [0,0.5]^2$ with $N_{\vx} = 100$ in each dimension; and the velocity domain is $(v_1,v_2) \in [-L_{\vv},L_{\vv}]^2$ with $L_{\vv} = 7.86$ and $N_{\vv}=32$. The Maxwell diffusive boundary condition is assumed at all boundaries. The wall quantities at $x_2 = 0.5$ are $\vu_{w} = (1,0)$, $T_{w} = 1$, while at all other boundaries we set $\vu_{w} = (0,0)$, $T_{w} = 1$. For the initial condition, we use a spatially homogeneous Maxwellian with $\rho_0 = 1$, $\vu_0= (1,1)$ and $T_0= 1$. 

\begin{figure}[!htbp]
\begin{center}
\includegraphics[width=2.4in]{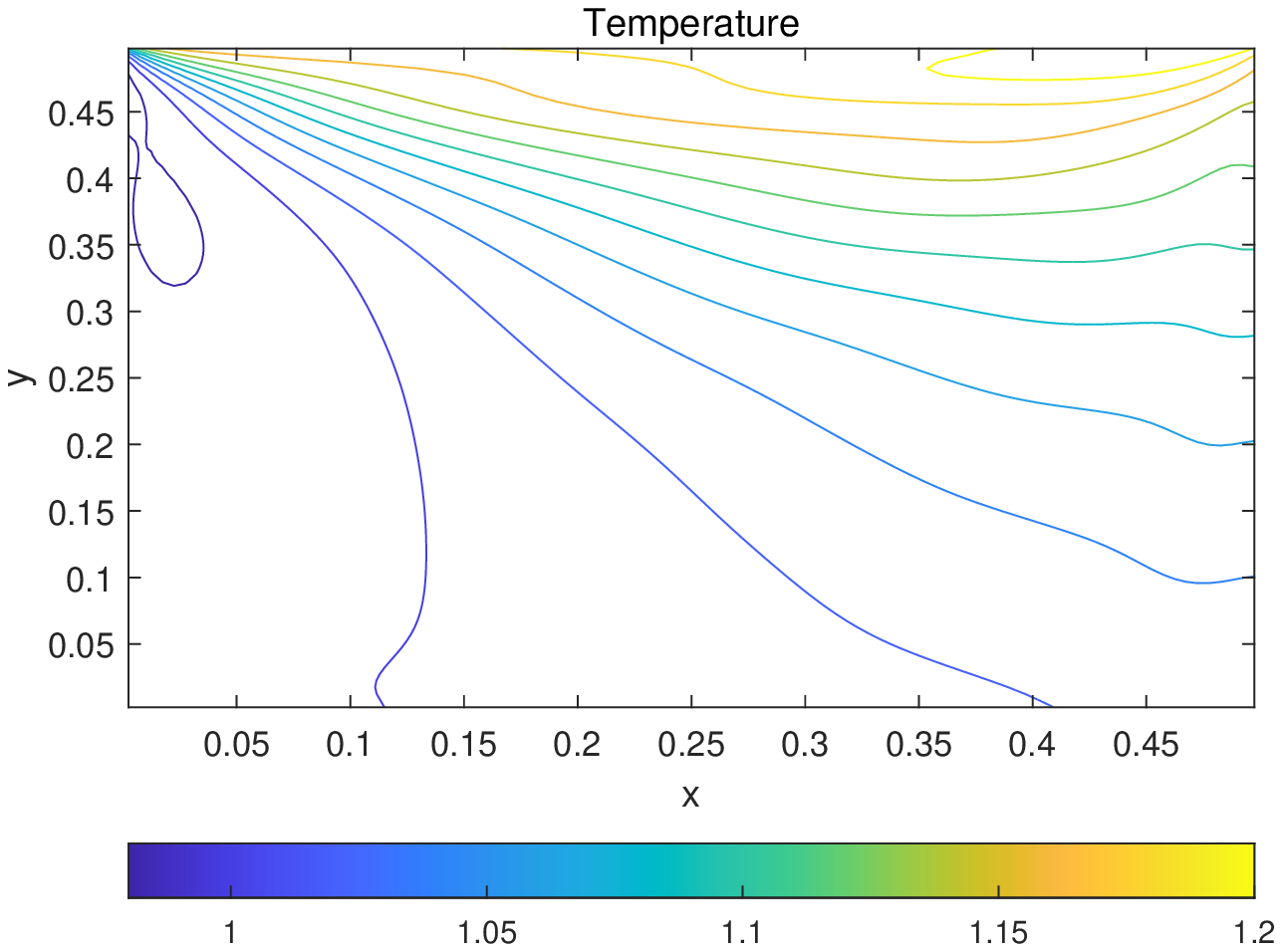}
\includegraphics[width=2.4in]{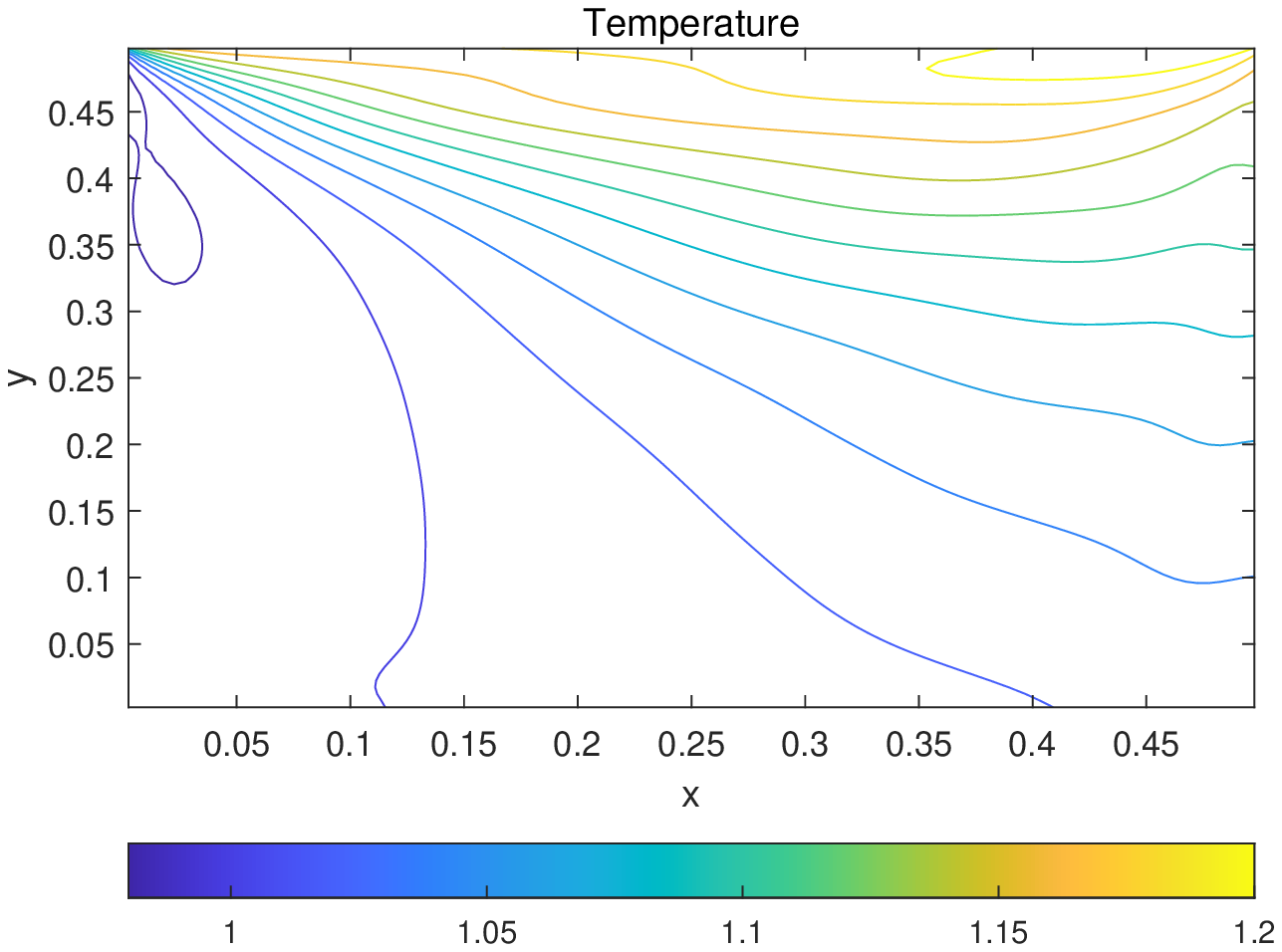}
\includegraphics[width=2.4in]{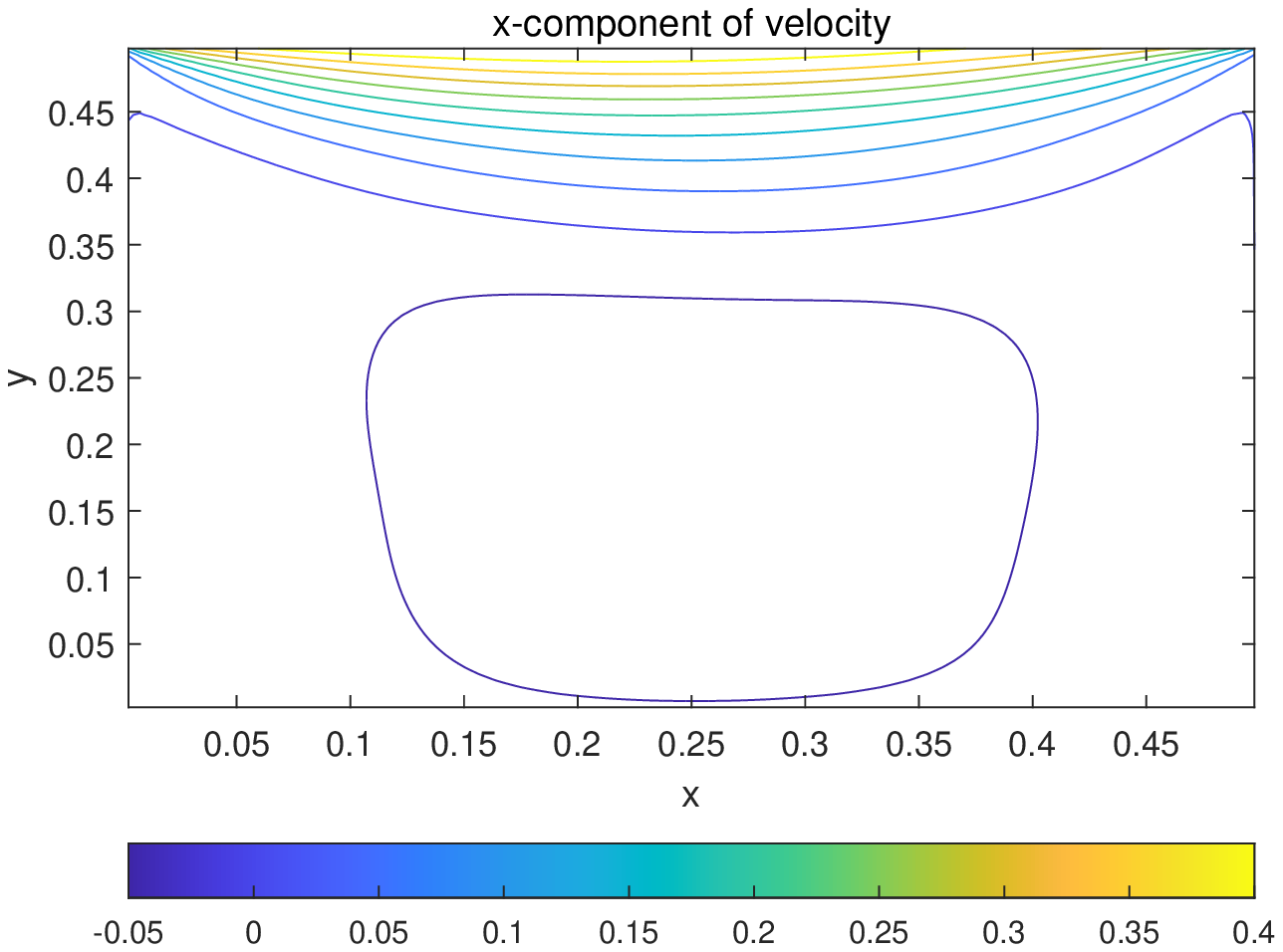}
\includegraphics[width=2.4in]{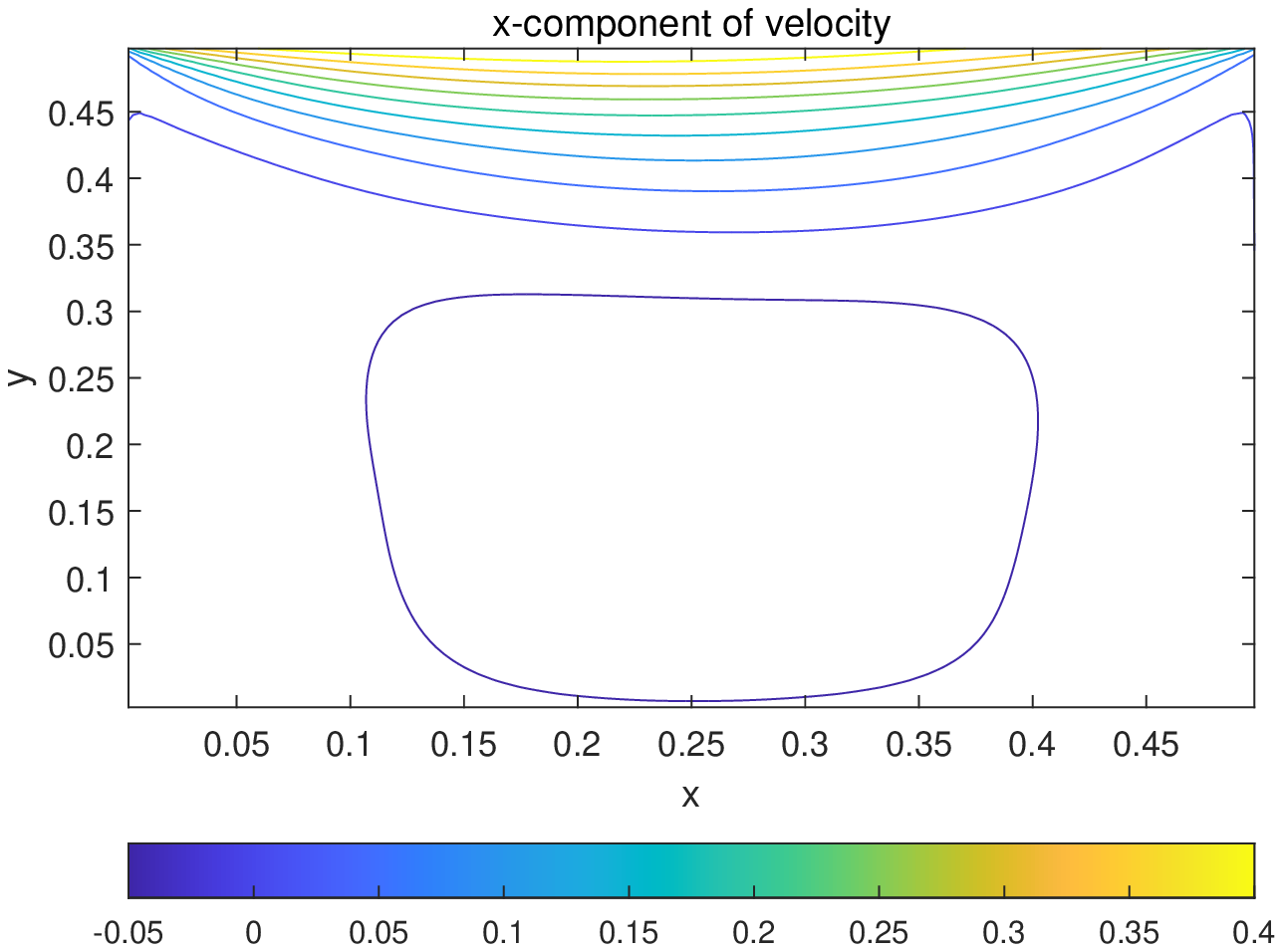}
\includegraphics[width=2.4in]{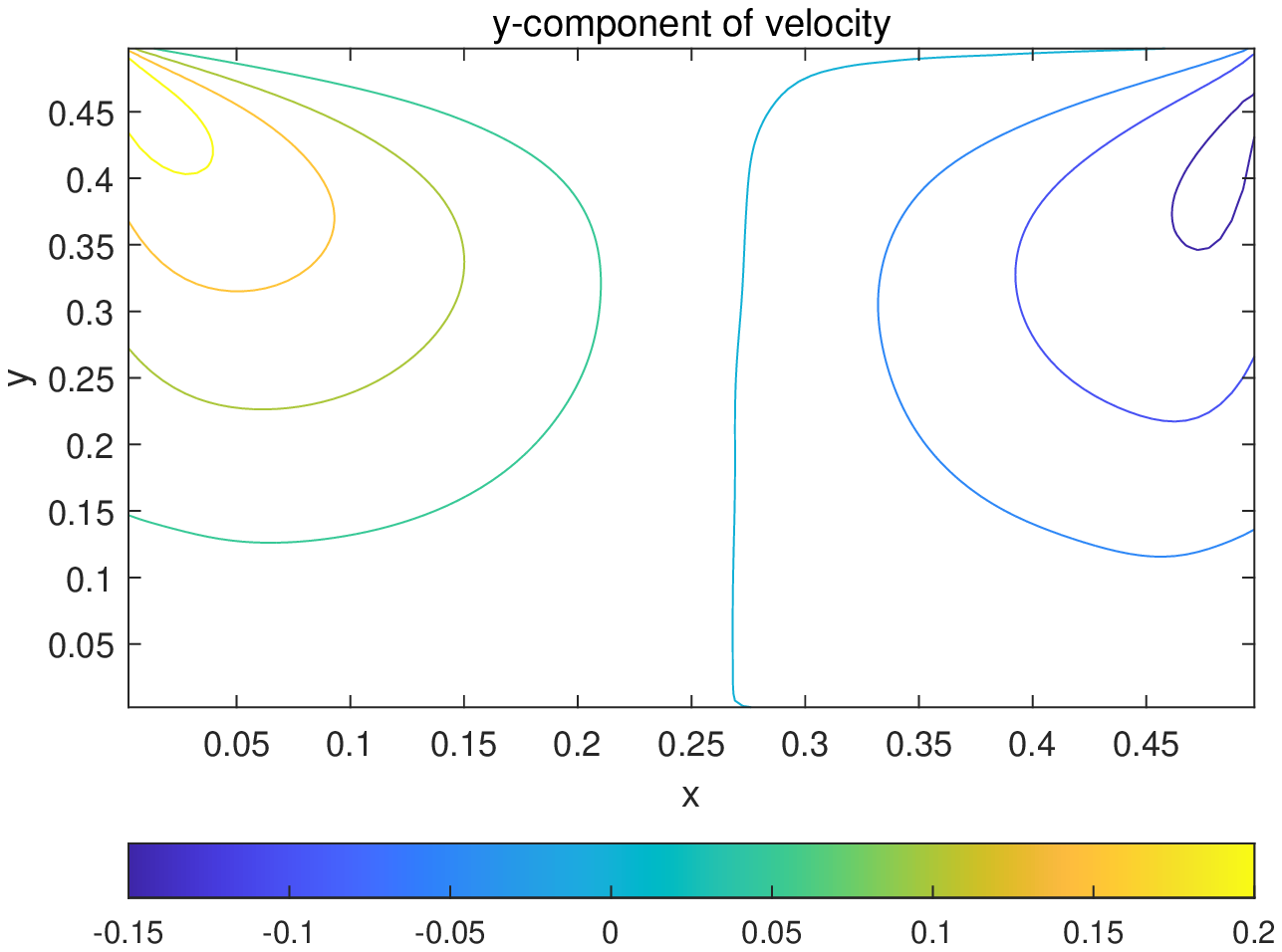}
\includegraphics[width=2.4in]{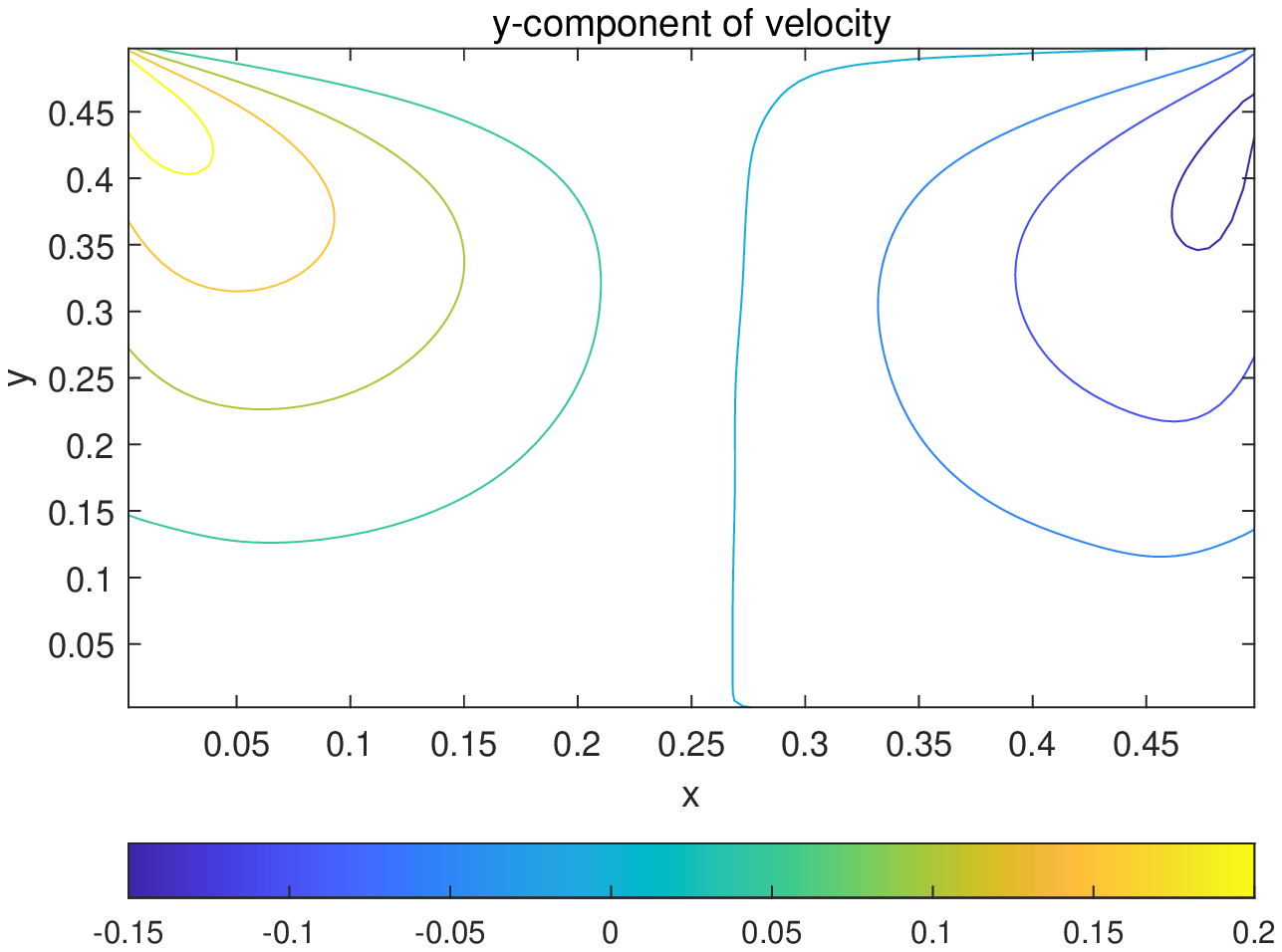}
\caption{Lid driven cavity flow. Top: temperature profile of the full grid method (Left) and low rank method (Right); Middle: $x_1$-component velocity of the full grid method (Left) and low rank method (Right); Bottom: $x_2$-component velocity of the full grid method (Left) and low rank method (Right). Convergence criterion is $\text{\bf res\_tol} = 2\times 10^{-7}$ for both methods.}
\label{ldcf mac}
\end{center}
\end{figure}

The convergence criterion is set as $\text{\bf res\_tol} = 2\times 10^{-7}$ for both the full grid method and the adaptive low rank method. For the full grid method, we need 29043 seconds to reach convergence, while for the low rank method, we only need 8323 seconds. We compare the temperature and velocity profile in \cref{ldcf mac} and a good match is obtained. 

\begin{figure}[!htbp]
\begin{center}
\includegraphics[width=2.4in]{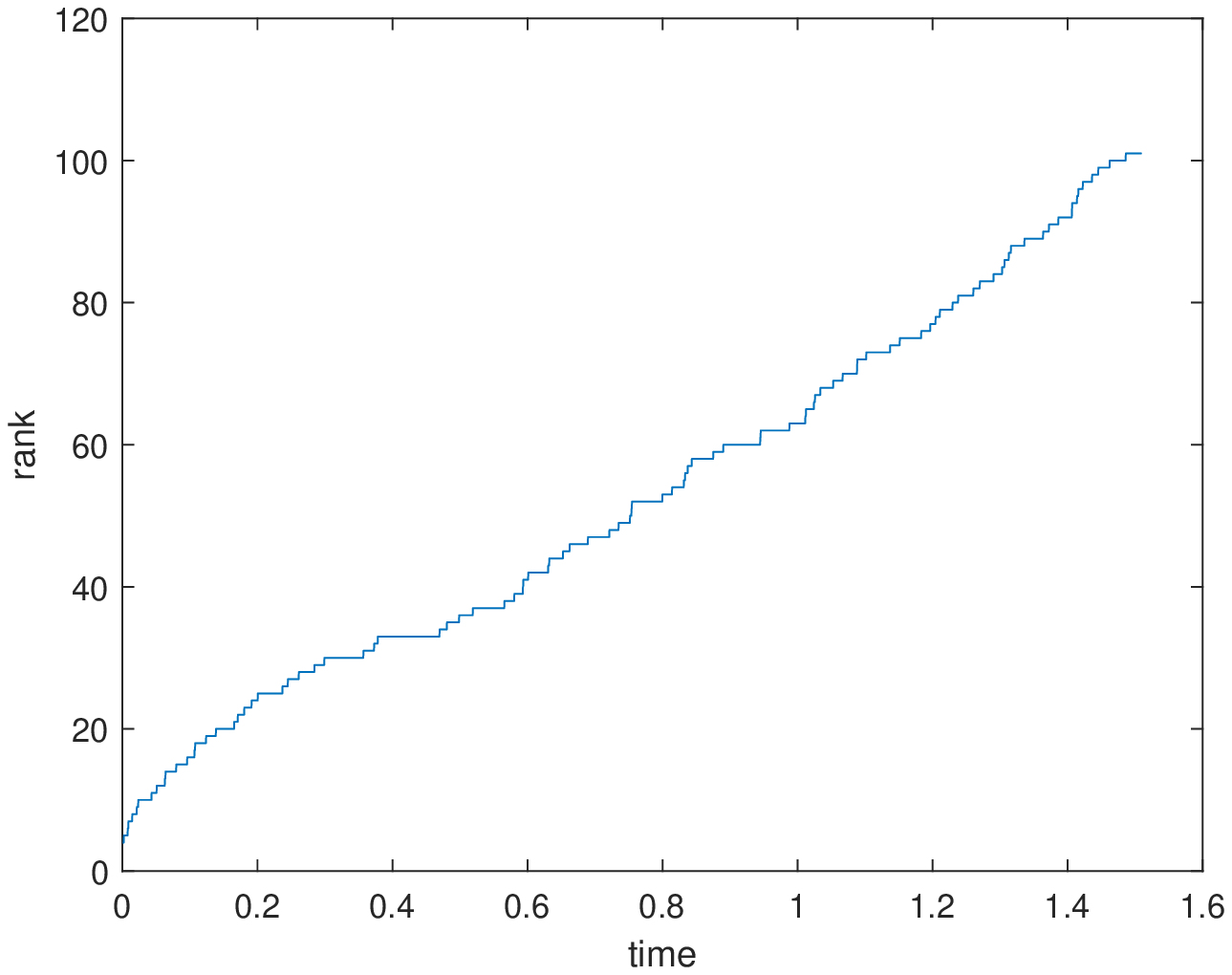}
\includegraphics[width=2.4in]{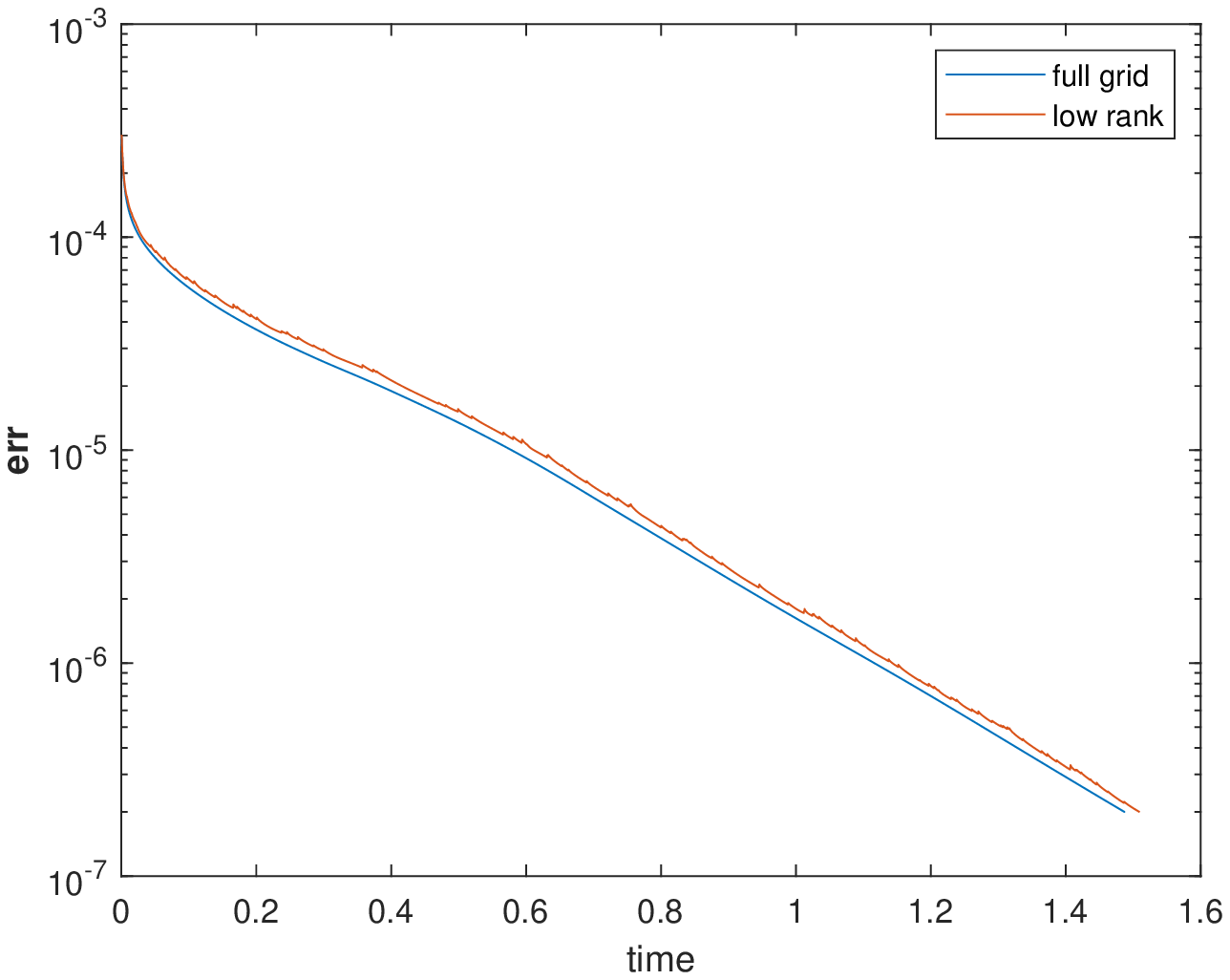}
\caption{Lid driven cavity flow. Left: rank evolution in the adaptive low rank method; Right: error decaying behaviors of the full grid method ($\text{err}_\text{\text{full tensor}}$) and low rank method ($\text{err}_\text{\text{low rank}}^{\text{ada}}$).}
\label{ldcf rank/err}
\end{center}
\end{figure}

From \cref{ldcf rank/err}, we can see that the rank in the adaptive low rank method is increasing with time and no stabilization is observed here, which implies this is an intrinsically high rank problem. Nevertheless, the error decay in the adaptive low rank method behaves similarly as in the full grid method (so our adaptive procedure does reasonable things in the actual simulation).


\subsection{Thermally driven cavity flow}

\begin{figure}[!htbp]
\begin{center}
\includegraphics[width=2.4in]{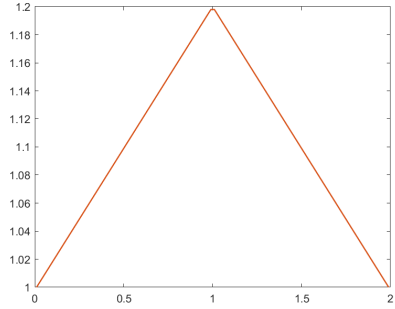}
\caption{Thermally driven cavity flow. Wall temperature profile at $x_2= 0$ and $x_2= 2$.}
\label{tdcf wall temp}
\end{center}
\end{figure}

We finally consider the 2D flow induced by thermal gradients. The spatial domain is rectangular $(x_1,x_2) \in [0,2]^2$ with $N_{\vx} = 100$ in each dimension; and the velocity domain is $(v_1,v_2) \in [-L_{\vv},L_{\vv}]^2$ with $L_{\vv} = 6.55$ and $N_{\vv} = 32$. The Maxwell diffusive boundary condition is assumed at all boundaries. We set the wall quantities at $x_2 = 0$, $x_2 = 2$ with $\vu_{w} = (0,0)$ and $T_{w}$ follows a linear function ranging from 1 to 1.2 as in \cref{tdcf wall temp}. At $x_1 = 0$, $x_1 = 2$, the wall quantities are set with $\vu_{w} = (0,0)$ and $T_{w} = 1$. For the initial condition, we use a spatially homogeneous Maxwellian with $\rho_0= 1$, $\vu_0= (0,0)$ and $T_0= 1$. 

\begin{figure}[!htbp]
\begin{center}
\includegraphics[width=2.4in]{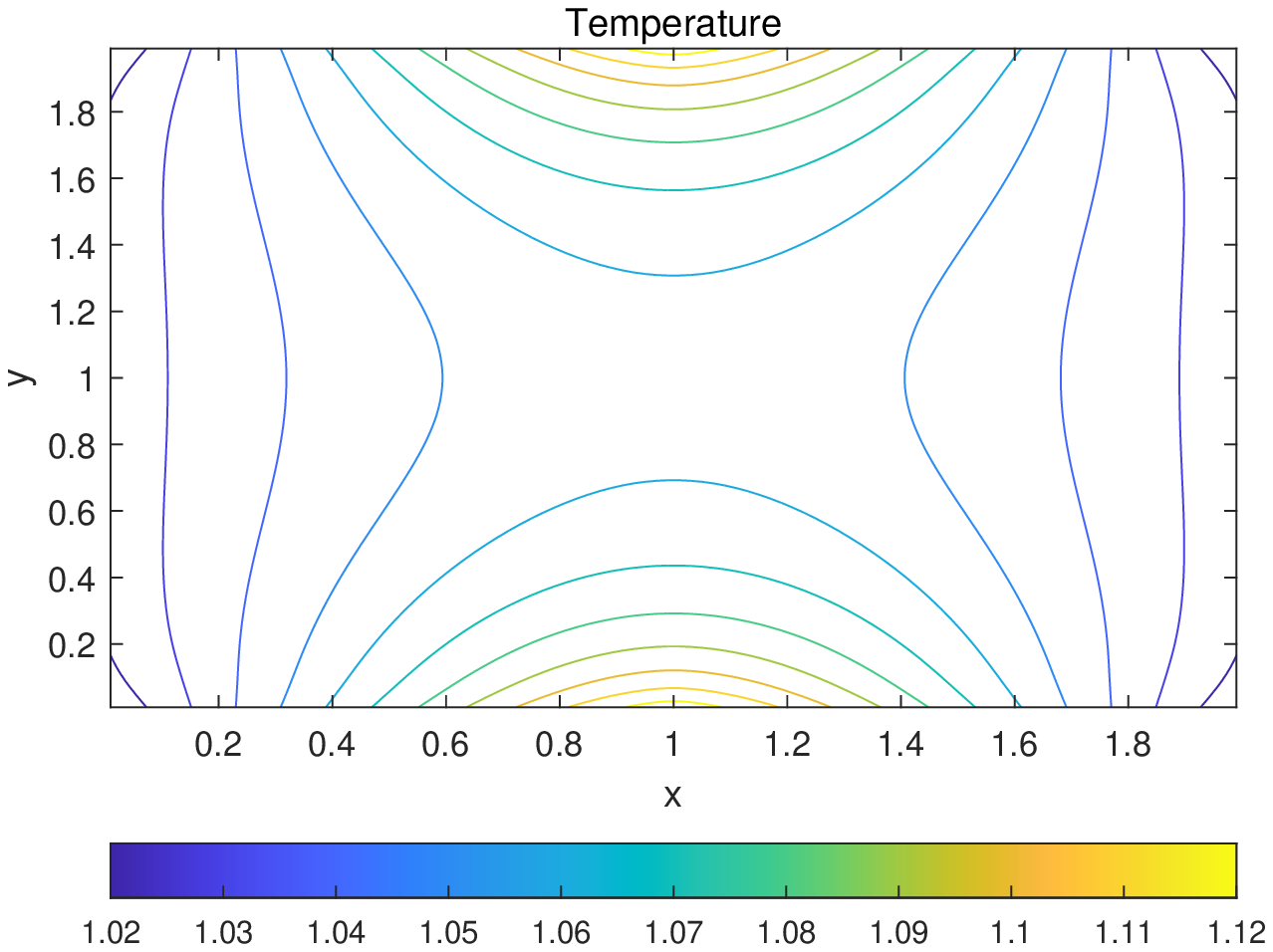}
\includegraphics[width=2.4in]{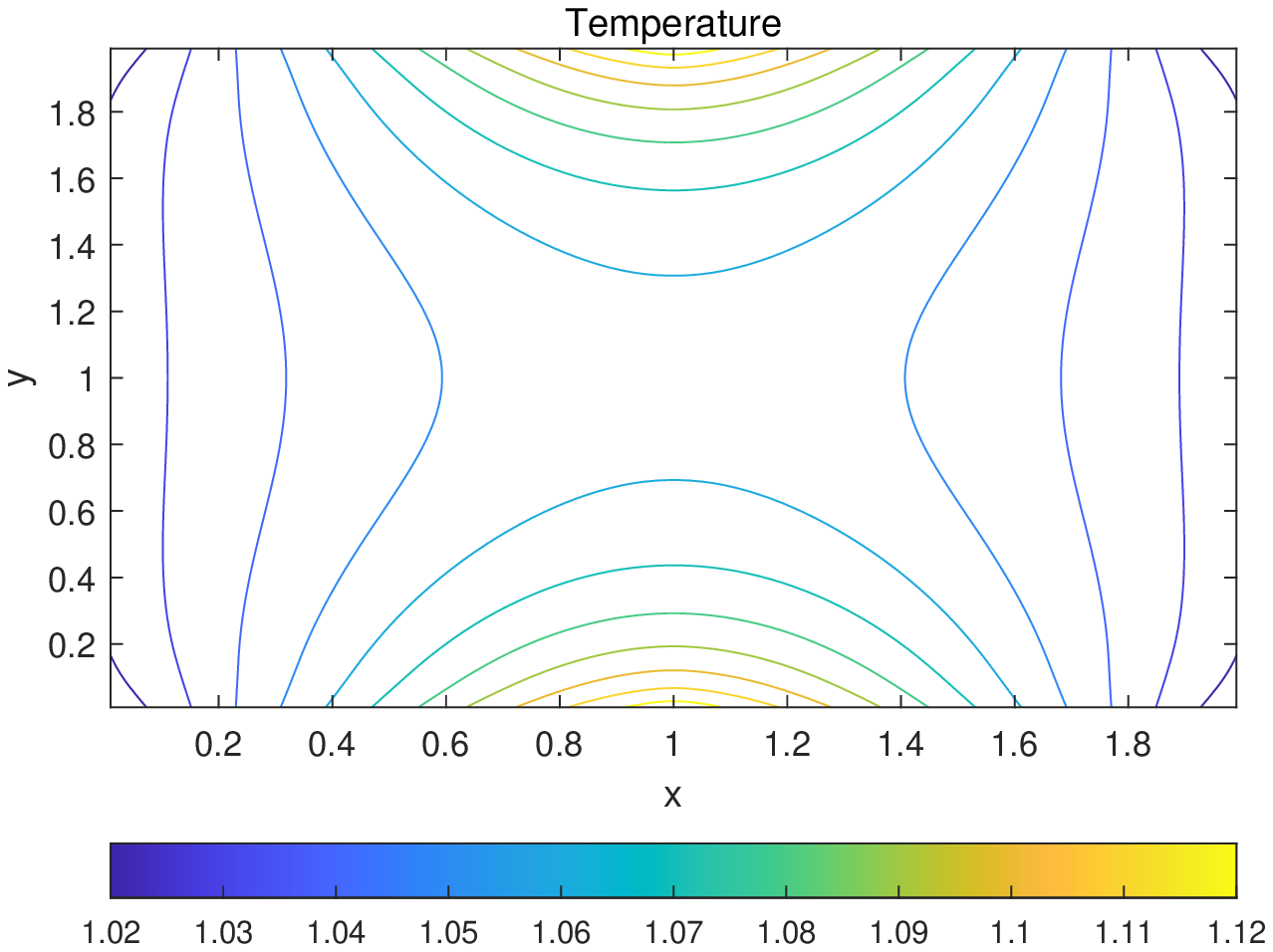}
\includegraphics[width=2.4in]{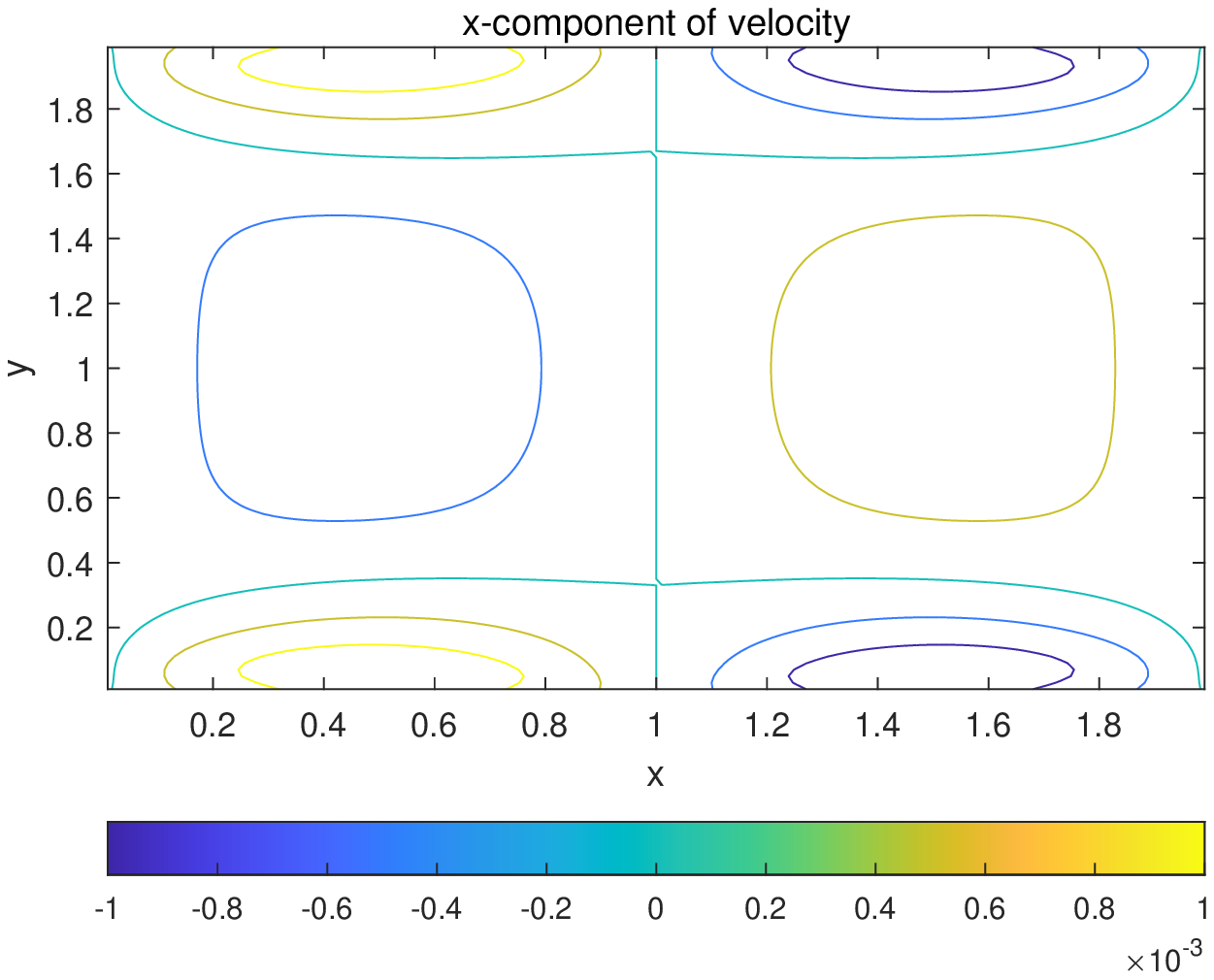}
\includegraphics[width=2.4in]{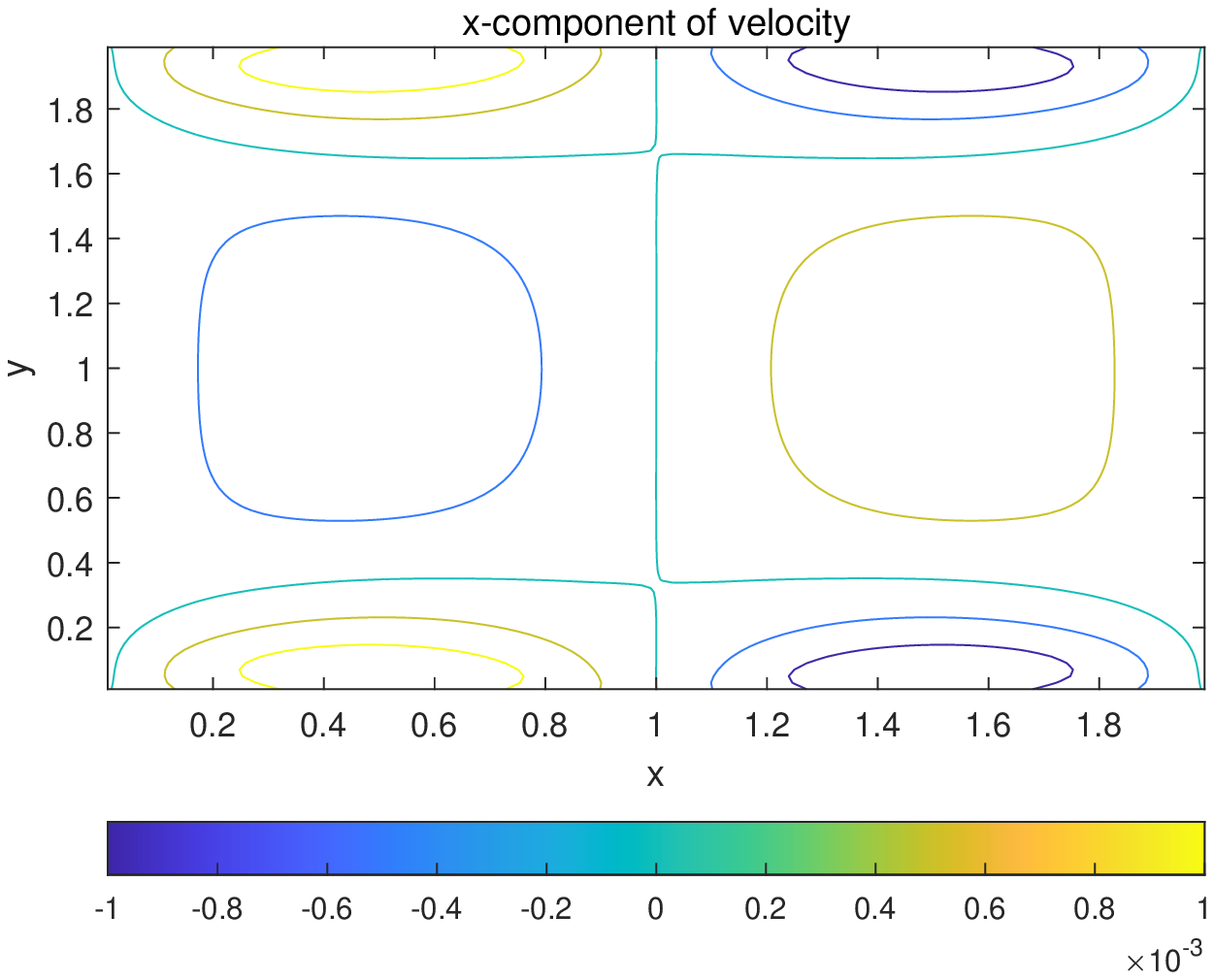}
\includegraphics[width=2.4in]{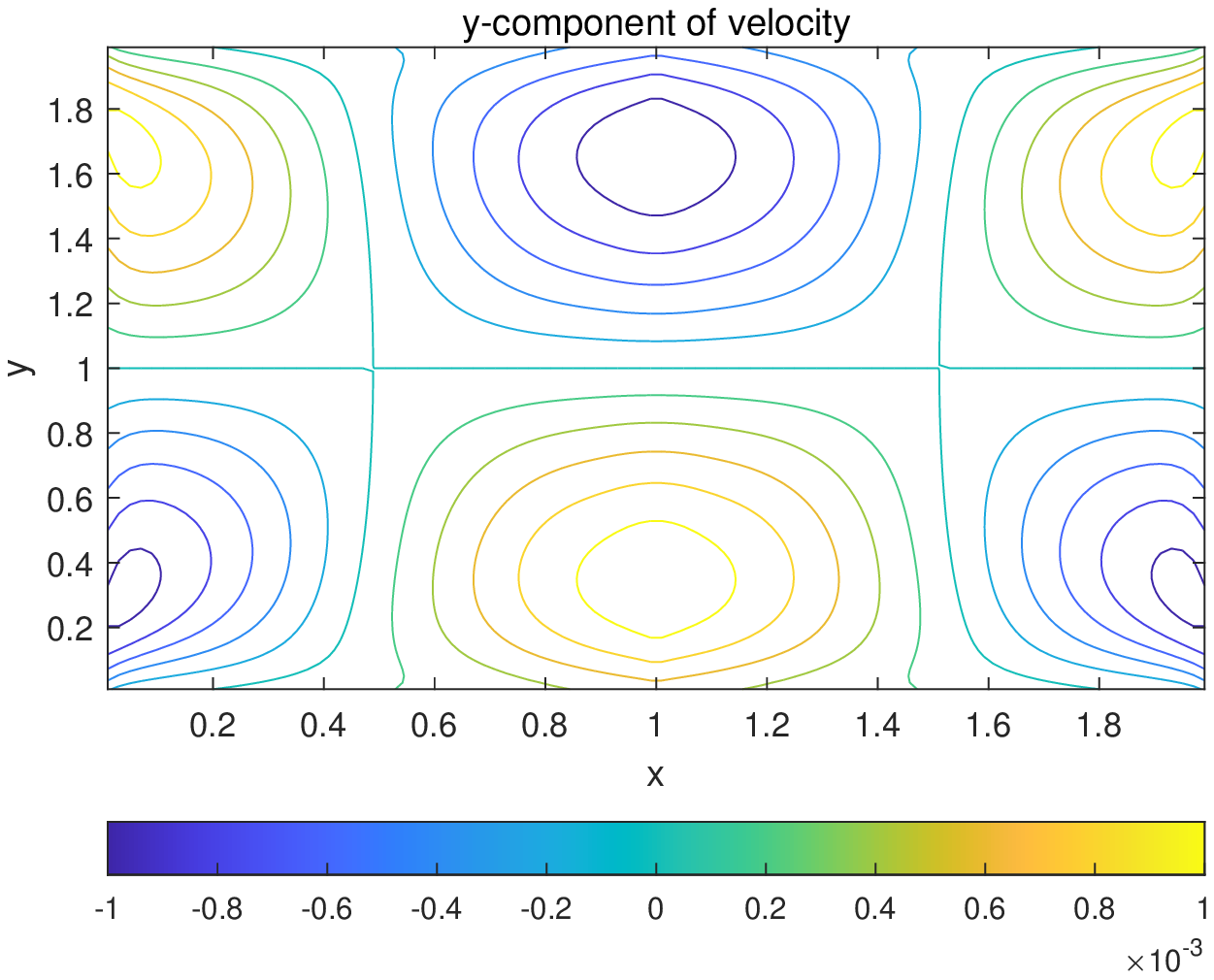}
\includegraphics[width=2.4in]{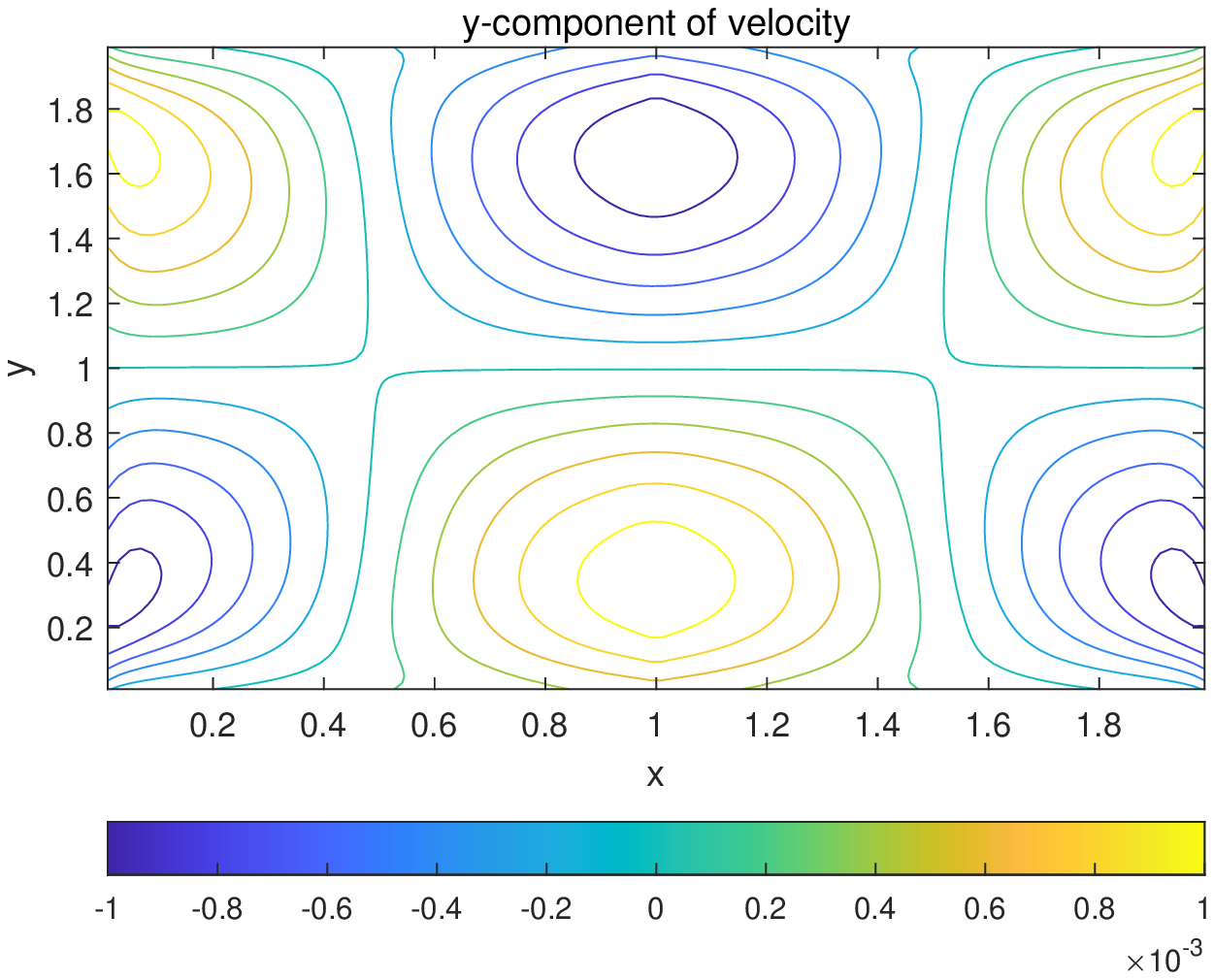}
\caption{Thermally driven cavity flow. Top: temperature profile of the full grid method (Left) and low rank method (Right); Middle: $x_1$-component velocity of the full grid method (Left) and low rank method (Right); Bottom: $x_2$-component velocity of the full grid method (Left) and low rank method (Right). Convergence criterion is $\text{\bf res\_tol} = 2\times 10^{-7}$ for both methods.}
\label{tdcf mac}
\end{center}
\end{figure}

The convergence criterion is set as $\text{\bf res\_tol} = 2\times 10^{-7}$ for both the full grid method and the adaptive low rank method. For the full grid method, we need 19011 seconds to reach convergence criterion, while for the low rank method, we only need 7112 seconds. We plot the temperature and velocity profile for both methods as in \cref{tdcf mac} where we can see a good match. 

\begin{figure}[!htbp]
\begin{center}
\includegraphics[width=2.4in]{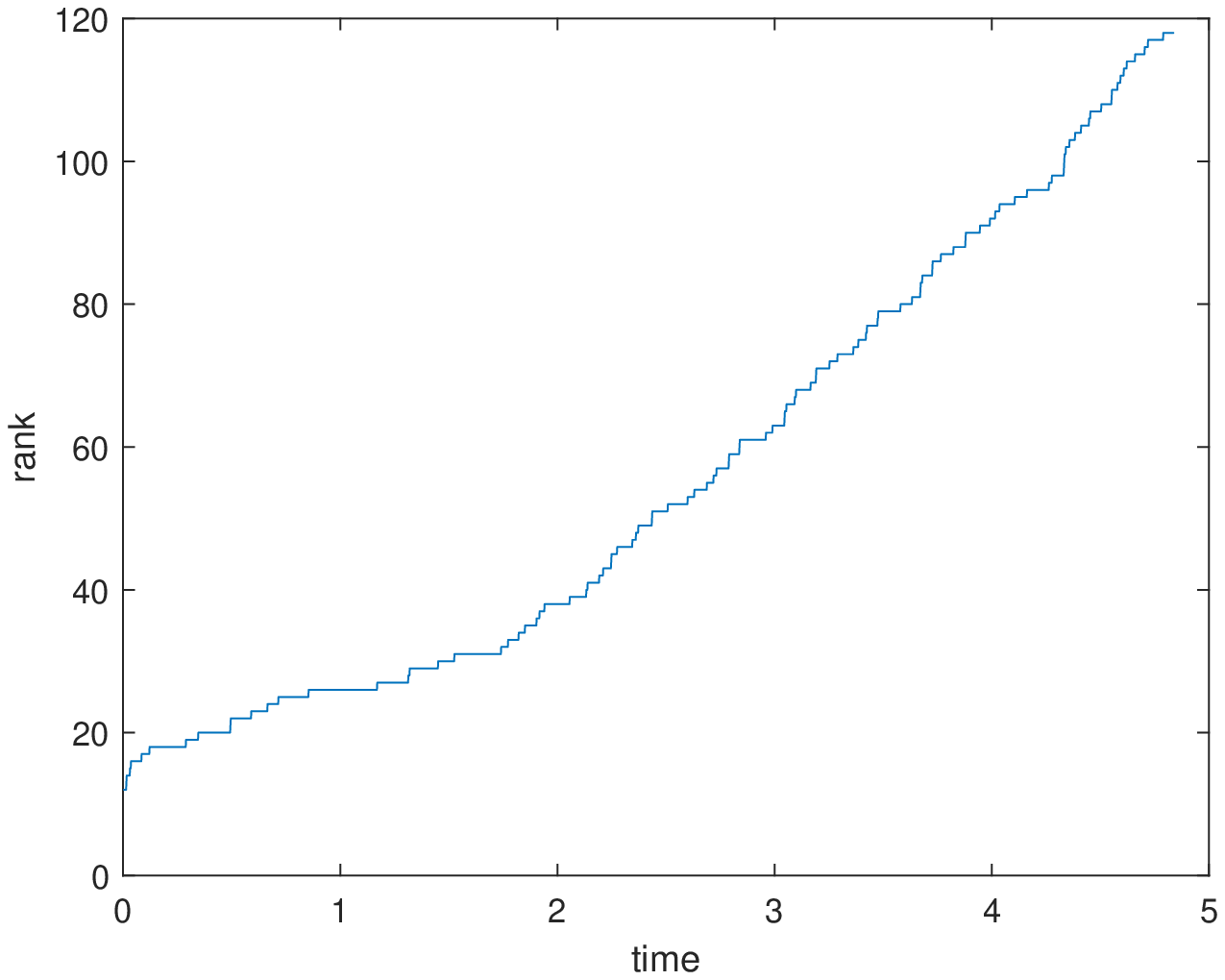}
\includegraphics[width=2.4in]{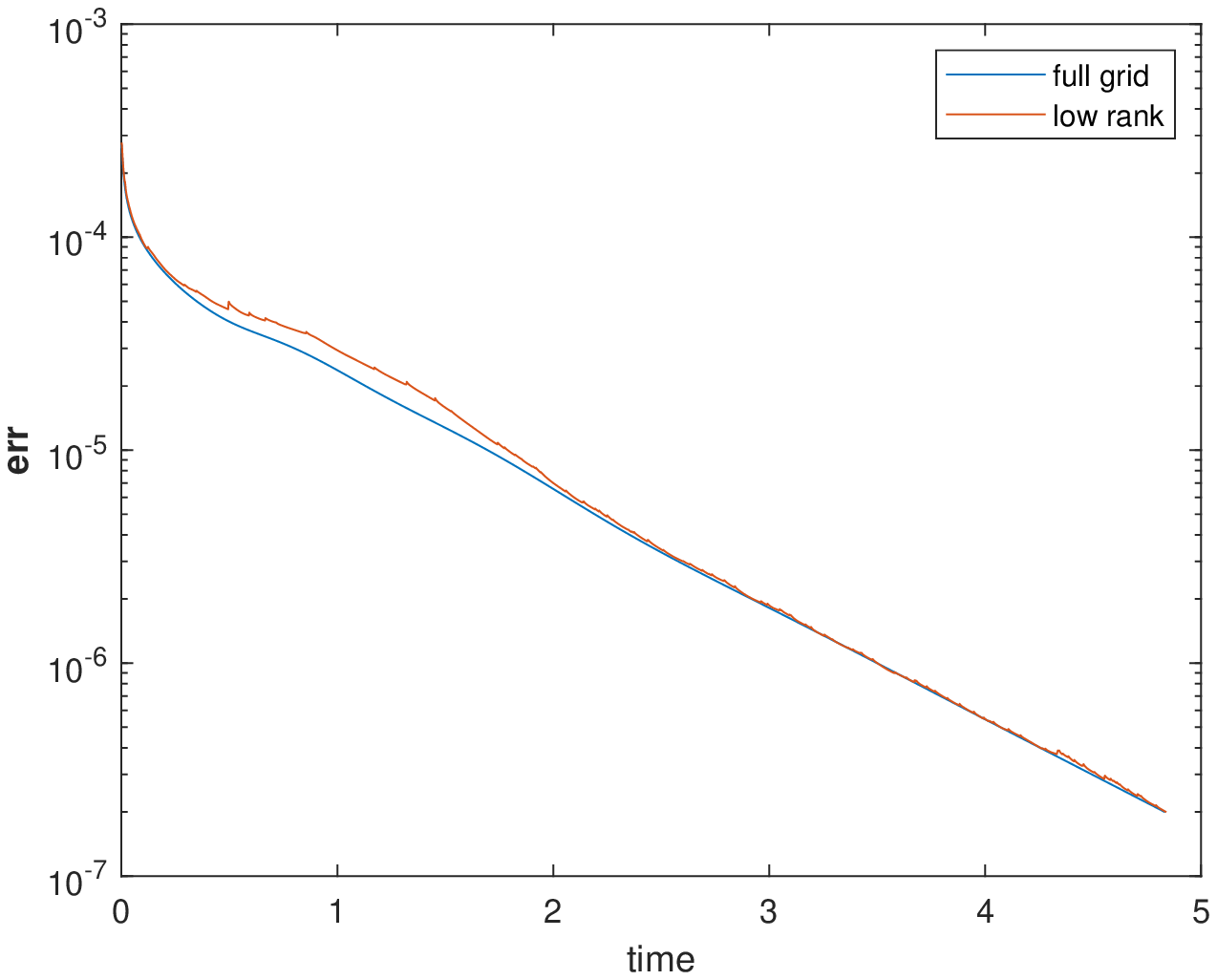}
\caption{Thermally driven cavity flow. Left: rank evolution in the adaptive low rank method; Right: error decaying behaviors of the full grid method ($\text{err}_\text{\text{full tensor}}$) and low rank method ($\text{err}_\text{\text{low rank}}^{\text{ada}}$).}
\label{tdcf rank/err}
\end{center}
\end{figure}

Similarly as in the previous test, we track the rank evolution in the adaptive low rank method and the error decay behavior of both methods in \cref{tdcf rank/err}. For this problem, the rank increases more rapidly, yet the low rank method can still produce reasonable solution more efficiently compared to the full grid method.


\section{Conclusions}
\label{sec:con}

We have introduced an adaptive dynamical low rank method for the nonlinear Boltzmann equation, concerning in particular the steady state computation. This method employs the fast Fourier spectral method (for the collision operator) and the dynamical low rank method to obtain computational efficiency. An adaptive strategy was introduced to incorporate the boundary information and control the computational rank by monitoring the residual error. A series of benchmark tests were performed to demonstrate the efficiency and accuracy of the proposed method in comparison to the full tensor grid method.

\section*{Data availability}

This manuscript has no associated data.

\bibliographystyle{plain}
\bibliography{normal_shock_ref.bib,hu_bibtex.bib}
\end{document}